\documentclass[aos,MSNbibl,dvips]{arximspdf}
\usepackage{graphicx}
\doi{10.1214/13-AOS1145} 
\volume{41}
\issue{5}
\pubyear{2013}
\firstpage{2324}
\lastpage{2358}

\makeatletter

\newcommand{\rright}{\right}
\newcommand{\lleft}{\left}
\newcommand{\eqref}[1]{(\ref{#1})}

\newtheorem{Thm}{Theorem}

\newproclaim{Defa}{Definition}
\newtheorem{Lem}[Thm]{Lemma}

\newproclaim{Def}{Definition}
\newproclaim{Con}{Condition}

\newproclaim{Post}{Postulate}
\newproclaim{Ex}{Example}

\newproclaim{Rem}{Remark}

\def\implies{\Longrightarrow}
\def\Exp{{\mathbb{E}}}

\newcommand{\ci}{\mathfrak{C}}
\newcommand{\oi}{\mathfrak{O}}

\newcommand{\Z}{\mathbb{Z}}

\def\independent{{\perp\!\!\!\!\perp}}
\makeatother

\begin{document}
\begin{frontmatter}

\title{Quantifying causal influences}
\runtitle{Quantifying causal influences}

\begin{aug}
\author[a]{\fnms{Dominik}~\snm{Janzing}\corref{}\thanksref{tt1}\ead[label=e1]{dominik.janzing@tuebingen.mpg.de}},
\author[b]{\fnms{David}~\snm{Balduzzi}\thanksref{tt1,tt2}\ead[label=e2]{david.balduzzi@inf.ethz.ch}},
\author[a]{\fnms{Moritz}~\snm{Grosse-Wentrup}\thanksref{tt1}\ead[label=e3]{moritz.grosse-wentrup@tuebingen.mpg.de}}
\and
\author[a]{\fnms{Bernhard}~\snm{Sch\"olkopf}\thanksref{tt1}\ead[label=e4]{bernhard.schoelkopf@tuebingen.mpg.de}}
\runauthor{Janzing, Balduzzi, Grosse-Wentrup and Sch\"olkopf}
\affiliation{Max Planck Institute for Intelligent Systems\thanksmark{tt1}
and ETH Z{\"u}rich\thanksmark{tt2}}
\address[a]{D. Janzing\\
M. Grosse-Wentrup\\
B. Sch\"olkopf\\
Max Planck Institute\\
\quad for Intelligent Systems\\
Spemannstr. 38\\
72076 T\"ubingen\\
Germany\\
\printead{e1}\\
\phantom{E-mail:\ }\printead*{e3}\\
\phantom{E-mail:\ }\printead*{e4}}

\address[b]{D. Balduzzi\\Max Planck Institute\\
\quad for Intelligent Systems\\
Spemannstr. 38\\
72076 T\"ubingen\\
Germany\\
and\\
ETH Z{\"u}rich\\
CAB F 63.1\\
Universit\"atstrasse 6\\
8092 Zurich\\
Switzerland\\
\printead{e2}}

\end{aug}

\received{\smonth{8} \syear{2012}}
\revised{\smonth{6} \syear{2013}}

%
\begin{abstract}
Many methods for causal inference generate directed acyclic graphs
(DAGs) that formalize causal relations between $n$ variables.
Given the joint distribution on all these variables, the DAG contains
all information about how intervening on one variable
changes the distribution of the other $n-1$ variables.
However, \textit{quantifying} the causal influence of one variable on
another one remains a nontrivial question.

Here we propose a set of natural, intuitive postulates that a measure
of causal strength should satisfy. We then introduce a communication
scenario, where edges in a DAG play the role of channels that can be
locally corrupted by interventions. Causal strength is then the
relative entropy distance between the old and the new distribution.

Many other measures of causal strength have been proposed, including
average causal effect, transfer entropy, directed information, and
information flow. We explain how they fail to satisfy the postulates on
simple DAGs of $\leq3$ nodes. Finally, we investigate the behavior of
our measure on time-series, supporting our claims with experiments on
simulated data.
\end{abstract}
%
%

\end{frontmatter}

\section{Introduction}\label{sec1}

Inferring causal relations is among the most important scientific goals
since causality, as opposed to mere statistical dependencies, provides
the basis for reasonable human decisions. During the past decade, it
has become popular to phrase causal relations in directed acyclic
graphs (DAGs) \cite{Pearl00} with random variables (formalizing
statistical quantities after repeated observations) as nodes and causal
influences as arrows. 

We briefly explain this formal setting. Here and throughout the paper,
we assume causal sufficiency, that is, there are no hidden variables
that influence more than one of the $n$ observed variables.
Let $G$ be a causal DAG with nodes $X_1,\ldots,X_n$ where
$X_i\rightarrow X_j$ means that $X_i$ influences $X_j$ ``directly'' in
the sense that
intervening on $X_i$ changes the distribution of $X_j$ even if all
other variables are held constant (also by interventions).
To simplify notation, we will mostly assume the $X_j$ to be discrete.
$P(x_1,\ldots,x_n)$ denotes the probability mass function of the joint
distribution $P(X_1,\ldots,X_n)$.
According to the Causal Markov Condition \cite{Spirtes,Pearl00}, which
we take for granted in this paper, every node $X_j$ is conditionally
independent of its nondescendants, given its parents
with respect to the causal DAG $G$.
If $PA_j$ denotes the set of parent variables of $X_j$ (i.e., its
direct causes) in $G$,
the joint probability thus factorizes \cite{Lauritzen1996} into
%
%
\begin{equation}
\label{fac} P(x_1,\ldots,x_n)=\prod
_{j=1}^n P(x_j|pa_j),
\end{equation}
where $pa_j$ denotes the values of $PA_j$.
By slightly abusing the notion of conditional probabilities, we assume
that $P(X_j|pa_j)$ is also defined for those $pa_j$ with $P(pa_j)=0$.
In other words, we know how the causal mechanisms act on potential
combinations of values of the parents that never occur.
Note that this assumption has implications because such causal
conditionals cannot be learned from observational data even if
the causal DAG is known.

Given this formalism, why define causal strength? After all, the DAG
together with the causal conditionals
contain the complete causal information: one can easily compute how
the joint distribution changes when an external intervention sets some
of the variables to specific values \cite{Pearl00}.
However,
describing causal relations in nature with a DAG always requires first
deciding how detailed the description should be. Depending on the
desired precision, one may want to account for some weak causal links
or not. Thus, an objective measure distinguishing weak arrows from
strong ones is required.

\subsection{Related work}

We discuss some definitions of causal strength that are either known or
just come up as straightforward ideas.

\emph{Average causal effect}: Following \cite{Pearl00}, $P(Y| {do}(X=x))$
denotes the distribution of $Y$ when $X$ is set to the value $x$ [it
will be introduced
more formally in equation~(\ref{do})].
Note that it only coincides with the usual conditional distribution
$P(Y| x)$ if the statistical dependence between $X$ and $Y$ is due to
a direct influence
of $X$ on $Y$, with no confounding common cause.
If all $X_i$ are binary variables, causal strength can then be quantified
by
the Average Causal Effect \cite{Holland,Pearl00}
\[
\operatorname{ACE}(X_i \rightarrow X_j):= P\bigl(X_j=1|
{do} (X_i=1)\bigr)- P\bigl(X_j=1| {do} (X_i=0)
\bigr).
\]
If a real-valued variable $X_j$ is affected by a binary variable $X_i$,
one considers
the shift of the mean of $X_j$ that is caused by switching $X_i$ from
$0$ to $1$.
Formally, one considers the difference \cite{Northcott}
\[
\Exp\bigl(X_j |{do} (X_i=1)\bigr)-\Exp
\bigl(X_j| {do} (X_i=0)\bigr).
\]
This measure only accounts for the linear aspect of an interaction
since it does not reflect whether $X_i$ changes higher order moments of
the distribution of $X_j$.

\emph{Analysis of Variance \textup{(}ANOVA\textup{)}}:
Let $X_i$ be caused by $X_1,\ldots,X_{i-1}$. The
variance of $X_i$ can formally be split into the average of the
variances of $X_i$, given $X_k$ with $k\leq i-1$, and the
variance of the expectations of $X_i$, given~$X_k$:
\begin{equation}\label{eq2}
\operatorname{Var} (X_i) = \Exp\bigl(\operatorname{Var} (X_i|X_k)
\bigr) + \operatorname{Var} \bigl(\Exp(X_i|X_k)\bigr).
\end{equation}
In the common scenario of drug testing experiments, for instance, the
first term in equation (\ref{eq2}) is given by the variability
of $X_i$ within a group of equal treatments (i.e., fixed $x_k$), while
the second one describes how much the means of $X_i$ vary between
different treatments.
It is tempting to say that the latter describes the part of the total
variation of $X_i$ that is \textit{caused by} the variation of $X_k$,
but this is conceptually wrong for nonlinear influences
and if there are statistical dependencies between $X_k$ and the other
parents of $X_i$ \cite{Lewontin,Northcott}.

For linear structural equations,
\[
X_i =\sum_{j<i} \alpha_{ij}
X_j +E_i \qquad\mbox{with $E_j$ being jointly
independent},
\]
and additionally assuming $X_k$ to be independent of the other parents
of $X_i$, the second term is given by $\operatorname{Var}(\alpha_{ik} X_k)$,
which indeed
describes the amount by which the variance of $X_i$ decreases when
$X_k$ is set to a fixed value by intervention.
In this sense,
%
%
\begin{equation}
\label{rdef} r_{ik}:=\frac{ \operatorname{Var} (\alpha_{ik} X_k) }{
\operatorname{Var} (X_i) }
\end{equation}
is indeed the fraction of the variance of $X_i$ that is \textit{caused by}
$X_k$. By rescaling all $X_j$ such that $\operatorname{Var}(X_j)=1$,
we have
$r_{ik}=\alpha^2_{ik}$. Then,
the square of the structure coefficients itself can be seen as a
simple measure for causal strength.

\emph{\textup{(}Conditional\textup{)} Mutual information}: The information of $X$ on $Y$
or vice versa is given by \cite{cover}
\[
I(X;Y):=\sum_{x,y} P(x,y)\log\frac{P(x,y)}{P(x)P(y)}.
\]
The information of $X$ on $Y$ or vice versa if $Z$ is given is defined
by \cite{cover}
%
%
\begin{equation}
\label{cmi} I(X;Y |Z):=\sum_{x,y,z} P(x,y,z)\log
\frac{P(x,y|z)}{P(x|z)P(y|z)}.
\end{equation}
There are situations where these expressions (with $Z$ describing some
background condition)
can indeed be interpreted as measuring the strength of the arrow
\mbox{$X\rightarrow Y$}. An essential part of this paper
describes the conditions where this makes sense and how to replace the
expressions with other information-theoretic ones
when it does not.

\emph{Granger causality/Transfer entropy/Directed information}:
Quantifying\break causal influence between time series [e.g. between
$(X_t)_{t\in\Z}$ and $(Y_t)_{t\in\Z}$] is special because one is
interested in quantifying the effect of all
$(X_t)$ on all $(Y_{t+s})$. If we represent the causal relations by a
DAG where every time instant defines a separate pair of variables, then
we ask for the strength of a \textit{set of arrows}. If $X_t$ and $Y_t$
are considered as instances of the variables $X,Y$, we leave the regime
of i.i.d. sampling.

Measuring the reduction of uncertainty in one variable after knowing
another is also a key idea in several related methods for quantifying
causal strength in time series. Granger causality in its original formulation
uses reduction of variance \cite{Granger1969}. Nonlinear
information-theoretic extensions in the same spirit are transfer
entropy~\cite{Schreiber} and directed information \cite{Massey}.
Both are
essentially based on conditional mutual information, where each
variable $X,Y,Z$ in (\ref{cmi}) is replaced with an appropriate set of
variables.



\emph{Information flow}: Since the above measures quantify dependencies
rather than
causality, several authors have defined causal strength by replacing
the observed probability distribution with distributions that arise
after interventions
(computed via the causal DAG).
\cite{AyKrakauer} defined Information Flow via an operation, ``source
exclusion'', which removes the influence of a variable in a network.
\cite{AyInfoFlow} defined a different notion of Information Flow
explicitly via Pearl's ${do}$-calculus.
Both measures are close to ours in spirit and in fact the version in
\cite{AyKrakauer} coincides with ours when quantifying the strength of
a single arrow. However, both do not satisfy our postulates.

\emph{Mediation analysis}:
\cite{Pearlindirect,Avin,Robins} explore how to separate the influence
of $X$ on $Y$ into parts that can be attributed to specific paths
by ``blocking'' other paths.
Consider, for instance, the case where $X$ influences $Y$ directly
and indirectly via $X\rightarrow Z\rightarrow Y$. To test its direct
influence, one changes $X$ from some ``reference'' value $x'$ to an
``active'' value $x$ while keeping the distribution of $Z$ that
either corresponds to the reference value $x'$ or
to the natural distribution $P(X)$.
A natural distinction between a reference state and an active state occurs,
for instance, in drug testing scenario where taking the drug
means switching from reference to active.
In contrast, our goal is not to study the impact
of one specific switching from $x'$ to $x$. Instead, we want
to construct a measure that quantifies the direct effect of the
variable $X$ on $Y$,
while treating all possible values of $X$ in the same way.
Nevertheless, there are interesting relation between
these approaches and ours that we briefly discuss at the end of
Section~\ref{subsecse}.

%



\section{Postulates for causal strength}
\label{secpost}

Let us first discuss the properties we expect a measure of causal
strength to have. The key idea is that causal strength is supposed to measure
the impact of an intervention that removes the respective arrows.
We present five properties that we consider reasonable.
Let $\ci_S$ denote the
strength of the arrows in set $S$. By slightly overloading notation, we
write $\ci_{X\rightarrow Y}$ instead of $\ci_{\{X\rightarrow Y\}}$.

\begin{longlist}[P0.]

\item[P0.]
\emph{Causal Markov condition}:
If $\ci_{S}=0$, then the joint distribution satisfies the Markov condition
with respect to the DAG $G_S$ obtained by removing the arrows in $S$.


\item[P1.]
\emph{Mutual information}:
If the true causal DAG reads $X\rightarrow Y$, then
\[
\ci_{X\rightarrow Y}=I(X;Y).
\]


\item[P2.]
\emph{Locality}:
The strength of $X \rightarrow Y$ only depends on (1) how $Y$ depends
on $X$ and its other parents, and (2) the joint distribution of all
parents of $Y$.
Formally, knowing $P(Y|PA_Y)$ and $P(PA_Y)$ is sufficient to compute
$\ci_{X\rightarrow Y}$. For strictly positive densities, this is
equivalent to knowing $P(Y,PA_Y)$.

\item[P3.]
\emph{Quantitative causal Markov condition}:
If there is an arrow from $X$ to $Y$, then the causal influence of $X$
on $Y$ is greater than or equal to
the conditional mutual information between $Y$ and $X$ given all the
other parents of
$Y$. Formally
\[
\ci_{X\rightarrow Y} \geq I\bigl(X;Y |PA_Y^X\bigr).
\]

\item[P4.]
\emph{Heredity}:
If the causal influence of a set of arrows is zero, then the causal
influence of all its subsets (in particular, individual arrows) is also zero.
\[
\mbox{If }S\subset T, \mbox{ then }\ci_T = 0 \implies
\ci_S=0.
\]
\end{longlist}

Note that we do not claim that \textit{every} reasonable measure of causal
strength should satisfy these postulates, but we now explain why we consider
them natural and show that the postulates make sense for simple DAGs.

P0:
If the purpose of our measure of causal strength is to quantify
relevance of arrows, then removing a set of arrows with zero strength
must make no difference. If, for instance, $\ci_{X\rightarrow Y}=0$,
removing $X\rightarrow Y$ should not yield a DAG that is ruled out by
the causal Markov condition.

We should emphasize that $\ci_S$ can be nonzero even if $S$ consists
of arrows each individually having zero strength.


P1: The mutual information actually measures the strength
of statistical dependencies. Since all these dependencies are generated
by the influence of $X$ on $Y$
(and not by a common cause or $Y$ influencing $X$),
it makes sense to measure causal strength by strength of dependencies.
Note that
mutual information $I(X;Y)=H(Y)-H(Y|X)$ also quantifies the
variability in $Y$ that is due to the variability in $X$, see also
Section~\ref{scontrol}.


\textit{Mutual information versus channel capacity.}
Given the premise that causal strength should be an information-like
quantity, a natural alternative to mutual information is the capacity
of the information channel $x \mapsto P(Y|x)$, that is, the maximum
over all
values of mutual information $I_{Q(X)}(X;Y)$ for all input
distributions $Q(X)$ of $X$ when keeping the conditional $P(Y|X)$.

While mutual information $I(X;Y)$ quantifies the observable
dependencies, channel capacity quantifies the strength of the strongest
dependencies
that can be generated using the information channel $P(Y|X)$. In this
sense, $I(X;Y)$ quantifies the \textit{factual} causal influence, while
channel capacity
measures the \textit{potential} influence.
Channel capacity also accounts for the impact of setting $x$ to values
that rarely or never occur in the observations.
However, this sensitivity regarding effects of rare inputs can
certainly be a problem for estimating the effect from sparse data. We
therefore prefer mutual information $I(X;Y)$ as it better assesses the
extent to which \emph{frequently observed changes} in $X$ influence $Y$.

P2:
Locality implies that we can ignore causes of $X$ when computing $\ci
_{X\rightarrow Y}$, unless they are at the same time direct causes of $Y$.
Likewise, other effects of $Y$ are irrelevant.
Moreover, it does not matter \textit{how} the dependencies between
the parents are generated (which parent influences which one or whether
they are effects of a common cause), we only need to know their joint
distribution with $X$.

\begin{figure}[t]

\includegraphics{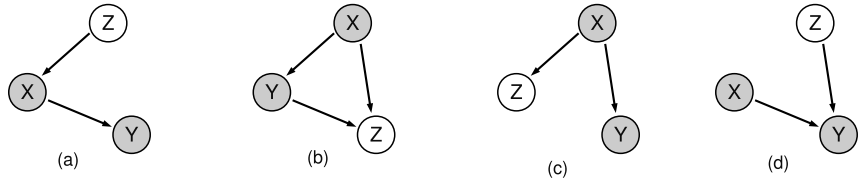}

\caption{DAGs for which the (conditional) mutual information is a
reasonable measure of causal strength:
For \textup{(a)} to \textup{(c)}, our postulates imply $\ci_{X\rightarrow Y} =I(X;Y)$. For
\textup{(d)} we will obtain
$\ci_{X\rightarrow Y}=I(X;Y |Z)$. The nodes $X$ and $Y$ are shaded
because they are source
and target of the arrow $X\rightarrow Y$, respectively.}\label{DAGsobv}
\end{figure}

Violations of locality have paradoxical implications. Assume, for example,
variable $Z$ would be relevant in DAG \ref{DAGsobv}(a). Then,
$\ci_{X\rightarrow Y}$ would depend on the mechanism that generates the
distribution of $X$,
while we are actually concerned with the information flowing from $X$
to $Y$ instead of that flowing \textit{to} $X$ from other nodes.
Likewise, [see DAGs \ref{DAGsobv}(b) and \ref{DAGsobv}(c)] it is
irrelevant whether $X$ and $Y$ have further effects.

%
%

P3: To justify the name of this postulate, observe that
the restriction of P0 to the single arrow case
$S=\{X\rightarrow Y\}$ is equivalent to
\[
\ci_{X\rightarrow Y}=0\quad\Longrightarrow\quad I\bigl(Y;X |PA_Y^X
\bigr)=0.
\]
To see this, we use the ordered Markov condition \cite{Pearl00},
Theorem~1.2.6, which is known to be equivalent to the Markov condition mentioned
in the \hyperref[sec1]{Introduction}. It states that every node is
conditionally
independent of
its predecessors (according to some ordering consistent with the DAG),
given its parents. If $PR_Y$ denotes the predecessors of $Y$
for some ordering that is consistent with $G$ and $G_S$,
the ordered Markov condition for $G_S$ holds iff
%
%
\begin{equation}
\label{omc} Y \independent PR_Y | PA_Y^X,
\end{equation}
since the conditions for all other nodes remain the same as in $G$.
Due to the semi-graphoid axioms (weak union and contraction rule \cite
{Pearl00}),
(\ref{omc}) is equivalent to
\[
Y \independent PR_Y \setminus\{X\} | PA_Y \wedge Y
\independent X | PA_Y^X.
\]
Since the condition on the left is guaranteed by the Markov condition
on $G$, the Markov condition on $G_S$ is equivalent to $I(Y;X |PA_Y^X)=0$.

In words, the arrow $X\rightarrow Y$ is the only reason for the
conditional dependence $I(Y;X |PA_Y^X)$ to be nonzero, hence it is
natural to postulate that its strength cannot be smaller than the
dependence that it generates. Section~\ref{secproperties} explains
why we should not postulate
equality.

P4: The postulate provides a compatibility condition: if
a set of arrows has zero causal influence, and so can be eliminated
without affecting the causal DAG, then the same should hold for all
subsets of that set. We refer to this as the heredity property by
analogy with matroid theory, where heredity implies that every subset
of an independent set is independent.



\section{Problems of known definitions}

Our definition of causal strength is presented in Section~\ref
{secdef}. This section discusses problems with alternate measures of
causal strength.

\subsection{ACE and ANOVA}

The first two measures are ruled out by P0. Consider a relation
between three binary variables $X,Y,Z$, where $Y=X\oplus Z$ with $X$
and $Z$ being unbiased and independent. Then changing $X$ has no
influence on the statistics of $Y$. Likewise, knowing $X$ does not
reduce the variance of $Y$.
To satisfy P0, we need modifications that account for the fact that we do observe
an influence of $X$ on $Y$ for each fixed value $z$ although this influence becomes invisible after
marginalizing over $Z$.


\subsection{Mutual information and conditional mutual information}
\label{subsecfailuremi}

It suffices to consider a few simple DAGs to illustrate why mutual
information and conditional mutual information are \emph{not} suitable
measures of causal strength in general.

\begin{figure}[t]

\includegraphics{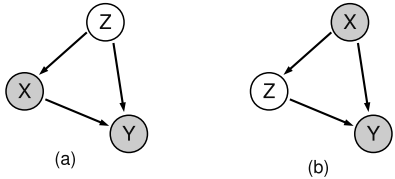}

\caption{DAGs for which finding a proper definition of $\ci
_{X\rightarrow Y}$ is challenging.}
\label{DAGs}
\end{figure}

\emph{Mutual information is not suitable in Figure~}\ref{DAGs}(a).
It is clear that $I(X;Y)$ is inappropriate because we can obtain
$I(X;Y)\neq0$ even when the arrow $X\rightarrow Y$ is missing, due to
the common cause $Z$. 

%

\emph{Conditional mutual information is not suitable for Figure}~\ref
{DAGs}(a).
%
Consider the limiting case\vadjust{\goodbreak}
where the direct influence $Z\rightarrow Y$ gets weaker until it almost
disappears ($P(y|x,z)\approx P(y|x)$).
Then the behavior of the system (observationally and interventionally)
is approximately described by the DAG~\ref{DAGsobv}(a). Using
$I(X;Y|Z)$ makes no sense in this scenario since, for example, $X$ may
be obtained from $Z$ by a simple copy operation, in which case
$I(X;Y|Z)=0$ necessarily, even when $X$ influences $Y$ strongly.


%
%
\begin{figure}[b]

\includegraphics{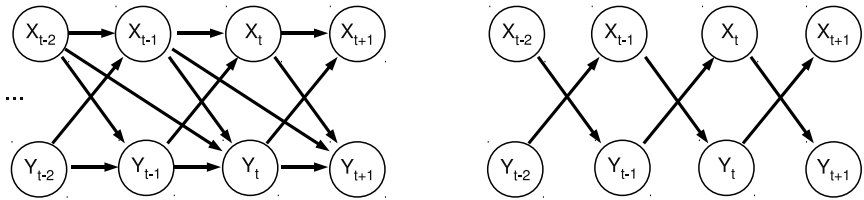}

\caption{Left: Typical causal DAG for two time series with mutual
causal influence.
The structure is acyclic because
instantaneous influences are excluded. Right: counter example in \cite
{AyInfoFlow}.
Transfer entropy vanishes if all arrows are copy operations although
the time series strongly influence each other.}\label{figTen}
\end{figure}

\subsection{Transfer entropy}
\label{subsecte}

Transfer entropy \cite{Schreiber} is intended to measure the influence
of one time-series on another one.
Let $(X_t,Y_t)_{t\in Z}$ be a bivariate stochastic process where $X_t$
influence some $Y_s$ with $s>t$, see Figure~\ref{figTen}, left.
Then transfer entropy is defined as the following conditional mutual
information:
\[
I(X_{(-\infty,t-1]} \rightarrow Y_{t} |Y_{(-\infty,t-1]}):=
I(X_{(-\infty,t-1]}; Y_{t} |Y_{(-\infty,t-1]}).
\]
It measures the amount of information the past of $X$ provides about
the present of $Y$ given the past of $Y$.
To quantify causal influence by conditional information relevance is also
in the spirit of Granger causality, where information is usually
understood in the sense of the amount of reduction of the linear
prediction error.

\emph{Transfer entropy is an unsatisfactory measure of causal strength}.
\cite{AyInfoFlow} pointed out that transfer entropy fails to quantify
causal influence for the following toy model:
Assume the information from $X_t$ is perfectly copied to $Y_{t+1}$ and
the information from $Y_t$ to $X_{t+1}$ (see Figure~\ref{figTen}, right).
Then the past of $Y$ is already sufficient to perfectly predict the
present value of $Y$ and the past of $X$ does not provide any further
information. Therefore, transfer entropy vanishes although both
variables heavily influence one another. If the copy operation is
noisy, transfer entropy is nonzero and thus seems more reasonable,
but the quantitative behavior is still wrong (as we will argue in
Example~\ref{experturbed}).

\begin{figure}

\includegraphics{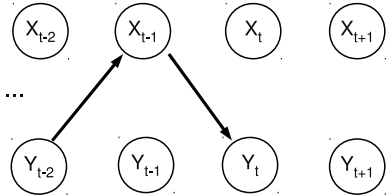}

\caption{Time series with only two causal arrows, where transfer
entropy fails satisfying our postulates.}\label{figTred}\vspace*{-3pt}
\end{figure}

\emph{Transfer entropy violates the postulates.}
Transfer entropy yields $0$ bits of causal influence in a situation
where common sense and P1 together with P2 require that
causal strength is $1$ bit (P2 reduces the DAG to one in which
P1 applies). Since our postulates refer to the strength of a
\textit{single} arrow while
transfer entropy is supposed to measure the strength of all arrows from
$X$ to $Y$, we reduce the DAG such that there is only one arrow from
$X$ to $Y$; see Figure~\ref{figTred}.
Then,
\begin{eqnarray*}
I(X_{(-\infty,t-1]} \rightarrow Y_{t} |Y_{(-\infty, t-1]})&=&
I(X_{(-\infty,t-1] };Y_{t} |Y_{(-\infty,t-1]})\\
&=&I(X_{t-1};
Y_{t} |Y_{t-2}).
\end{eqnarray*}
The causal structure coincides with DAG \ref{DAGsobv}(a) by setting
$Y_{t-2}\equiv Z$, $X_{t-1} \equiv X$, and $Y_{t}\equiv Y$.
With these replacements, transfer entropy
yields $I(X;Y |Z)=0$ bits instead of $I(X;Y)=1$ bit, as required by
P1 and P2.

%
%

Note that the same problem occurs if causal strength between time
series is quantified by directed information \cite{Massey} because this
measure also conditions on the entire past of $Y$.


\subsection{Information flow}
\label{siflow}

Note that \cite{AyInfoFlow} and \cite{AyKrakauer} introduce two
different quantities, both called ``information flow.'' We consider
them in turn.

After arguing that transfer entropy does not properly capture the
strength of the impact of interventions,
\cite{AyInfoFlow} proposes to define causal strength using Pearl's
${do}$ calculus \cite{Pearl00}. Given a causal directed acyclic graph
$G$, Pearl computes the joint distribution obtained if variable $X_j$
is forcibly set to the value $x_j$ as
%
%
\begin{equation}
\label{do} P \bigl(x_1,\ldots,x_n | {do}
\bigl({x}'_j\bigr) \bigr):=\prod_{i\neq j}
P(x_i|pa_i)\cdot\delta_{x_j, {x}'_j}.
\end{equation}
Intuitively, the intervention on $X_j$ removes the dependence of $X_j$
on its parents and therefore replaces\vadjust{\goodbreak} $P(x_j|pa_j)$ with
the kronecker symbol.
Likewise, one can define interventions on several nodes by replacing
all conditionals with kronecker symbols.

Given three sets of nodes $X_A$, $X_B$ and $X_C$ in a directed acyclic
graph $G$, information flow is defined by
\begin{eqnarray*}
&&I \bigl(X_A\rightarrow X_B |{do} ({X}_C)
\bigr)
\\
&&\qquad:= \sum_{x_{\vphantom{x'_{a} }C},x_{\vphantom{x'_{a} }A},x_{\vphantom{x'_{a} }B}}P(x_C) P
\bigl(x_A|{do} ({x}_C) \bigr) P \bigl(x_B
|{do}({x}_A, {x}_C) \bigr)\\
&&\hspace*{38pt}\quad\qquad{}\times\log\frac{P (x_B |{do}({x}_A,{x}_C) )}{\sum_{x_A'}P (x_A'
|{do}({x}_C) )P (x_B |{do}({x}_A', {x}_C) )}.
\end{eqnarray*}
To better understand this expression, we first consider the case where
the set $X_C$ is empty.
Then we obtain
\[
I(X_A \rightarrow X_B):= \sum
_{x_{\vphantom{x'_{a} } A},x_{\vphantom{x'_{a} } B}} P(x_A) P\bigl(x_B|{do}
(x_A)\bigr) \log\frac{P(x_B|{do} (x_A))}{\sum_{x_A'} P(x'_A)P(x_B|{do} (x'_A))},
\]
which measures the mutual information between $X_A$ and $X_B$ obtained
when the information channel $x_A \mapsto P(X_B| {do} (x_A))$ is used
with the input
distribution $P(X_A)$.


\emph{Information flow, as defined in \cite{AyInfoFlow}, is an
unsatisfactory measure of causal strength.}
To quantify $X\rightarrow Y$ in DAGs \ref{DAGs}(a) and \ref{DAGs}(b)
using information flow, we may either choose $I(X\rightarrow Y)$ or
$I(X\rightarrow Y |{do}(Z))$.
Both choices are inconsistent with our postulates and intuitive expectations.

Start with $I(X\rightarrow Y)$ and DAG~\ref{DAGs}(a).
Let $X,Y,Z$ be binary with $Y:=X\oplus Z$ an XOR. Let $Z$ be an
unbiased coin toss and $X$ obtained from $Z$ by a faulty copy operation
with two-sided symmetric error.
One easily checks that $I(X\rightarrow Y)$ is zero in the limit of
error probability $1/2$ (making $X$ and $Y$ independent).
Nevertheless, dropping the arrow $X\rightarrow Y$ violates the Markov
condition, contradicting P0.
For error rate close to $1/2$, we still violate P3 because
$I(Y;X |Z)$ is close to $1$, while $I(X\rightarrow Y)$ is close to zero.
A similar argument applies to DAG \ref{DAGs}(b).

Now consider $I(X\rightarrow Y |{do}(Z))$. Note that it yields
different results for DAGs~\ref{DAGs}(a) and \ref{DAGs}(b)
when the joint distribution is the same, contradicting P2. This
is because $P(x|{do}(z))=P(x| z)$ for \ref{DAGs}(a), while
$P(x|{do}(z))=P(x)$ for \ref{DAGs}(b).
In other words, $I(X\rightarrow Y |{do}(Z))$ depends on the causal
relation \emph{between} the two causes $X$ and $Z$, rather than only on
the relation between causes and effects.

Apart from being inconsistent with our postulate, it is unsatisfactory
that $I(X\rightarrow Y |{do}(Z))$ tends to zero for the example above if the
error rate of copying $X$ from $Z$ in DAG \ref{DAGs}(a) tends to zero
(conditioned on setting $Z$ to some value, the information passed from
$X$ to $Y$ is zero because $X$ attains a fixed value, too).\vadjust{\goodbreak} In this
limit, $Y$ is always zero. Clearly, however, link $X\rightarrow Y$ is
important for explaining the behavior of the XOR: without the link,
the gate would not output ``zero'' for both $Z=0$ and $Z=1$.

\begin{figure}

\includegraphics{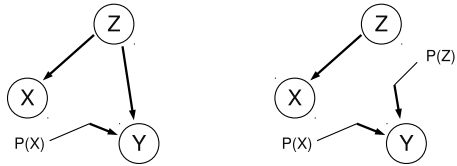}

\caption{Left: deletion of the arrow $X\rightarrow Y$. The conditional
$P(Y|X,Z)$ is fed with an independent copy of $X$, distributed with $P(X)$.
The resulting distribution reads $P_{X\rightarrow Y}(x,y,z)=P(x,z)\sum
_{x'} P(y|z,x')P(x')$.
Right: deletion of both incoming arrows. The conditional $P(Y|X,Z)$ is
then fed with the product distribution $P(X)P(Z)$ instead of the joint
$P(X,Z)$ as in \cite{AyKrakauer}, since the latter would require
communication between the open ends. We obtain $P_{X\rightarrow
Y,Z\rightarrow Y}(x,y,z)=\sum_{x',z'}P(x,z)P(y|x',z')P(x')P(z')$. Feeding
with independent inputs is particularly relevant for the following example:
let $X$ and $Z$ be binary with $X=Z$ and $Y=X\oplus Z$. Then, the
cutting had no impact
if we would keep the dependences.} \label{deletion}
\end{figure}

\emph{Information flow, as defined in \cite{AyKrakauer}, is
unsatisfactory as a measure of causal strength for sets of edges.}
Since this measure is close to ours,
we will explain (see caption of Figure~\ref{deletion}) the difference when
introducing ours and show that P4 fails without our modification.

%
%

\section{Defining the strength of causal arrows}
\label{secdef}

\subsection{Definition in terms of conditional probabilities}
\label{subsecdefcon}

This section proposes a way to quantify the causal influence of a set
of arrows that yields satisfactory answers in all the cases discussed above.
Our measure is motivated by a scenario where nodes represent different
parties communicating with each other via channels. Hence, we think of
arrows as physical channels that propagate information between distant
points in space, for example, wires that connect electronic devices.
Each such wire connects the output of a device with the input of
another one. For the intuitive ideas below, it is also important that
the wire connecting
$X_i$ and $X_j$ physically contains full information about $X_i$ [which
may be more than the information that is required to explain the output
behavior $P(X_j|PA_j)$].
We then think of the strength of arrow $X_i \rightarrow X_j$ as the
impact of corrupting it, that is, the impact of cutting the wire.
To get a well-defined ``post-cutting'' distribution we have to say what
to do with the open end corresponding to $X_j$, because it needs to be
fed with some input.
It is natural to feed it probabilistically with inputs $x_i$ according
to $P(X_i)$ because this is the only distribution of $X_i$ that is
locally observable
[feeding it with some conditional distribution $P(X_i|\cdot\cdot)$
assumes that
the one cutting the edge has access to other nodes---and not only the
physical state of the channel]. Note that this notion of cutting edges\vadjust{\goodbreak}
coincides with the ``source exclusion'' defined in \cite{AyKrakauer}
if only
one edge is cut.
However, we define the deletion of a \textit{set} of arrows by feeding all
open ends with the product of the corresponding marginal distributions, while
\cite{AyKrakauer} keeps the dependencies between the open ends and
removes the dependencies between open ends and the other variables.
Our post-cutting distribution can be thought of as arising from a scenario
where each channel is cut by an independent attacker, who
tries to blur the attack by feeding her open end with
$P(X_i)$ (which is the only distribution she can see), while
\cite{AyKrakauer} requires communicating attackers who agree on feeding
their open ends with the observed joint distribution.


Lemma~\ref{lemmarkovian} and Remark~\ref{remmarkovian} below provide
a more mathematical argument for the product distribution.
Figure~\ref{deletion} visualizes the deletion of one edge (left) and
two edges (right).


We now define the ``post-cutting'' distribution formally:\vspace*{-3pt}

%
\begin{Def}[(Removing causal arrows)]\label{dcisa}
Let $G$ be a causal DAG and $P$ be Markovian with respect to $G$. Let
$S\subset G$ be a set of arrows. Set $PA_j^S$ as the set of those
parents $X_i$ of $X_j$ for which $(i,j)\in S$ and $PA_j^{\bar{S}}$
those for which $(i,j)\notin S$. Set
%
%
\begin{equation}
P_S\bigl(x_j|pa^{\bar{S}}_j\bigr):=
\sum_{pa_j^S} P\bigl(x_j|pa^{\bar{S}}_j,pa_j^S
\bigr) P_{\prod}\bigl(pa_j^S\bigr), \label{epsj}
\end{equation}
where $P_{\prod} (pa_j^S)$ denotes for a given $j$ the product of
marginal distributions of all variables in $PA_j^{S}$.
Define a new joint distribution, the \emph{interventional
distribution}\footnote{Note that this intervention differs from the
kind of interventions considered by \cite{Pearl00}, where
variables are set to specific values. Here we intervene on the arrows,
the ``information channels,'' and not on the nodes.}
%
%
\begin{equation}
P_S(x_1,\ldots,x_n):=\prod
_j P_S\bigl(x_j|pa^{\bar{S}}_j
\bigr). \label{eps}\vspace*{-3pt}
\end{equation}
\end{Def}
%
See Figure~\ref{deletion}, left, for a simple example with cutting only
one edge.
Equation~(\ref{eps}) formalizes the fact that each open end of the wires
is independently fed with the corresponding marginal distribution, see also
Figure~\ref{deletion}, right.
Information flow in the sense of \cite{AyKrakauer} is obtained when the
product distribution $P_\Pi(pa_j^{\bar{S}})$
in (\ref{epsj}) is replaced with the joint distribution $P(pa_j^{\bar{S}})$.

The modified joint distribution $P_S$ can be considered as generated by
the reduced DAG:\vspace*{-3pt}
%
%
\begin{Lem}[(Markovian)]\label{lemmarkovian}
The interventional distribution $P_S$ is Markovian with respect to the
graph $G_S$ obtained from $G$ by removing the edges in $S$.\vspace*{-3pt}
\end{Lem}

\begin{pf}
By construction, $P_S$ factorizes according to $G_S$ in the sense of~(\ref{fac}).\vadjust{\goodbreak}
\end{pf}

%
\begin{Rem}
\label{remmarkovian}
Markovianity is violated if the dependencies between open ends are
kept. Consider, for instance, the DAG $X\rightarrow Y \rightarrow Z$.
Cutting both edges yields
\[
P_S(x,y,z)=P(x) \sum_{x'} P
\bigl(y|x'\bigr)P\bigl(x'\bigr) \sum
_{y'}P\bigl(z|y'\bigr)P\bigl(y'
\bigr)=P(x)P(y)P(z),
\]
which is obviously Markovian with respect to the DAG without arrows.
Feeding the ``open ends'' with $P(x',y')$ instead obtains
\[
\tilde{P}_S(x,y,z)=P(x)\sum_{x'y'} P
\bigl(y|x'\bigr)P\bigl(z|y'\bigr) P
\bigl(x',y'\bigr),
\]
which induces dependencies between $Y$ and $Z$, although we have
claimed to have removed all links between the three variables.
\end{Rem}

%
%
\begin{Def}[(Causal influence of a set of arrows)]\label{dciset}
The causal influence of the arrows in $S$ is given by the
Kullback--Leibler divergence
%
%
\begin{equation}
\ci_S(P):= D(P\|P_S). \label{ers}
\end{equation}
If $S=\{X_k\rightarrow X_l\}$ is a single edge we write
$\ci_{k\rightarrow l}$ instead of $\ci_{X_k \rightarrow X_l}$.
\end{Def}

%
\begin{Rem}[(Observing versus intervening)]
\label{observeintervene}
Note that $P_S$ could easily be confused with a different distribution
obtained when
the open ends are fed with conditional distributions rather than
marginal distributions.
As an illustrative example, consider DAG \ref{DAGs}(a) and define
$\tilde{P}_{X\rightarrow Y} (X,Y,Z)$ as
\[
\tilde{P}_{X\rightarrow Y}(x,y,z):=P(x,z)P(y|z) =P(x,z)\sum
_{x'} P\bigl(y|x'\bigr) P\bigl(x'|z
\bigr),
\]
and recall that replacing $P(x'|z)$ with $P(x')$ in the right most
expression yields $P_{X\rightarrow Y}$.
We call $\tilde{P}_{X\rightarrow Y} $ the ``partially observed
distribution.'' It is the distribution obtained by ignoring
the influence of $X$ on $Y$: $\tilde{P}_{X\rightarrow Y}$ is computed
according to (\ref{fac}), but uses a DAG where $X\rightarrow Y$ is missing.
The difference between ``ignoring'' and ``cutting'' the edge is
important for the following reason.
By a known rephrasing of mutual information as relative entropy \cite
{cover} we obtain
%
%
\begin{equation}
\label{parti} D(P\|\tilde{P}_{X\rightarrow Y})=\sum_{x,y,z}
P(x,y,z) \log\frac
{P(y|z,x)}{P(y|z)} = I(X;Y |Z),
\end{equation}
which, as we have already discussed, is \textit{not} a satisfactory
measure of causal strength. On the other hand, we have
%
%
\begin{eqnarray}
\label{rewr} \ci_{X\rightarrow Y}&=&D(P\|P_{X\rightarrow Y})=D \bigl
[ P(Y|Z,X)\|
P_{X\rightarrow Y} (Y|Z,X) \bigr]
\\
\label{rewr2}&=& D \bigl[P(Y|Z,X)\|P_{X\rightarrow
Y}(Y|Z) \bigr]
\nonumber
\\[-8pt]
\\[-8pt]
\nonumber
& = &\sum
_{x,y,z} P(x,y,z) \log\frac{P(y|z,x)}{\sum_{x'} P(y|z,x')P(x')}.
\end{eqnarray}
Comparing the second expressions in (\ref{rewr2}) and (\ref{parti})
shows again that the difference between ignoring and cutting is due to the
difference between $P(y|z)$ and $\sum_{x'} P(y|z,x')P(x')$.
\end{Rem}
The following scenario provides a better intuition for the rightmost
expression in (\ref{rewr2}).

%
\begin{Ex}[(Redistributing a vaccine)]
Consider the task of quantifying the effectiveness of a vaccine. Let
$X$ indicate whether a patient decides to get vaccinated or not and $Y$
whether the patient becomes infected. Further assume that the vaccine's
effectiveness is strongly confounded by age $Z$ because the vaccination
often fails for elderly people. At the same time, elderly people
request the vaccine more often because they are more afraid of infection.
Ignoring other confounders, the DAG in Figure~\ref{DAGs}(a) visualizes
the causal structure.

Deleting the edge $X\rightarrow Y$
corresponds to an experiment where the vaccine is randomly assigned to
patients regardless of their intent and age (while keeping the total
fraction of patients vaccinated constant).
Then $P_{X\rightarrow Y} (y|z,x)=P_{X\rightarrow Y}(y|z)= \sum_{x'}
P(y|z,x')P(x')$
represents the conditional probability of infection, given age, \emph
{when vaccines are distributed randomly}.
$\ci_{X\rightarrow Y}$ quantifies the difference to $P(y|z,x)$, which
is the conditional probability of infection, given age and intention
\emph{when patients act on their intentions}.
It thus measures the impact of destroying the coupling between the
intention to get the vaccine and getting it via randomized redistribution.
\end{Ex}

\subsection{Definition via structural equations}
\label{subsecse}

The definition above uses the conditional density $P(x_j|pa_j)$.
Estimating a conditional density from empirical data
requires huge samples or strong assumptions---particularly for
continuous variables. Fortunately, however, structural equations
(also called functional models \cite{Pearl00})
allow more direct estimation
of causal strength without referring to the conditional distribution.

%
\begin{Def}[(Structural equation)]
A structure equation is a model that explains the joint distribution
$P(X_1,\ldots,X_n)$ by
a deterministic dependence
\[
X_j=f_j(PA_j,E_j),
\]
where the variables $E_j$ are jointly independent unobserved noise
variables. Note that functions $f_j$ that correspond to parentless
variables can be chosen to be the identity, that is, $X_j=E_j$.
\end{Def}

Suppose that we are given a
causal inference method that directly infers the structural equations
(e.g., \cite{Hoyer,UAIidentifiability}) in the sense that
it outputs $n$-tuples $(e^i_1,\ldots,e^i_n)$ with $i=1,\ldots,m$ (with
$m$ denoting the sample size)
as well as the functions $f_j$ from the observed
$n$-tuples $(x^i_1,\ldots,x^i_n)$.

%
\begin{Def}[(Removing a causal arrow in a structural equation)]
Deletion of the arrow $X_k\rightarrow X_l$ is modeled by \textup{(i)}
introducing an i.i.d. copy $X'_k$ of $X_k$ and \textup{(ii)} subsuming the
new random variable $X'_k$ into the noise term of $f_l$. The result is
a new set of structural equations:
%
%
\begin{eqnarray}
\label{eseint} x_j & = &f_j (pa_j,e_j
) \qquad\mbox{if }j\neq l \mbox{ and }
\nonumber
\\[-8pt]
\\[-8pt]
\nonumber
x_l & =& f_l \bigl(pa_l\setminus
\{x_k\},\bigl(x'_k,e_l\bigr)
\bigr),
\end{eqnarray}
where we have omitted the superscript $i$ to simplify notation.
\end{Def}

%
\begin{Rem}
To measure the causal influence of a set of arrows, we apply the same
procedure after first introducing jointly independent i.i.d. copies of
all variables at the tails of deleted arrows.
\end{Rem}

%
\begin{Rem}
The change introduced by the deletion only affects
$X_l$ and its descendants, the virtual sample thus keeps all $x_j$ with
$j<l$. Moreover,
we can ignore all variables $X_j$ with $j>l$ due to Lemma~\ref{lemdown}.
\end{Rem}

Note that $x_k'$ must be chosen to be independent of all $x_j$ with
$j\leq k$, but, by virtue of the structural equations, not independent
of $x_l$ and its descendants. The new structural equations thus
generate $n$-tuples of ``virtual'' observations
$x^S_1,\ldots,x^S_n$ from the input
\[
\bigl(e_1,\ldots,\bigl(x'_k,e_l
\bigr),\ldots,e_n\bigr).
\]
We show below that $n$-tuples generated this way indeed follow the
distribution $P_S(X_1,\ldots,X_n)$.
We can therefore estimate causal influence via any method that
estimates relative entropy using
the observed samples $x_1,\ldots,x_n$ and the virtual ones $\tilde
{x}_1,\ldots,\tilde{x}_n$.
To illustrate the above scheme, we consider the case where $Z$ and $X$
are causes of $Y$ and we want to delete the edge $X\rightarrow Y$.
The case where $Y$ has more than $2$ parents follows easily.

%
\begin{Ex}[(Two parents)]
The following table corresponds to the observed variables $X,Z,Y$, as
well as the unobserved noise $E^Y$ which
we assumed to be estimated together with learning the structural
equations:
%
%
\begin{equation}
\label{eseed} \lleft(\matrix{ Z & X & E^Y & Y\vspace*{2pt}
\cr
\hline z_1 & x_1 & e_1^Y &
f_Y\bigl(z_1,x_1,e^Y_1
\bigr)\vspace*{2pt}
\cr
z_2 & x_2 & e_2^Y
& f_Y\bigl(z_2,x_2,e^Y_2
\bigr)\vspace*{2pt}
\cr
\vdots& & & \vdots\vspace*{2pt}
\cr
z_m &
x_m & e_m^Y & f_Y
\bigl(z_m,x_m,e^Y_m\bigr) }
\rright).
\end{equation}

To simulate the deletion of $X\rightarrow Y$ we first
generate a list of virtual observations for $Y$ after generating
samples from an i.i.d. copy $X'$ of $X$:
%
%
\begin{equation}
\label{eseedint} \lleft(\matrix{ Z & X & X' & E^Y
& Y\vspace*{2pt}
\cr
\hline z_1 & x_1 &
x_1' & e^Y_1 &
f_Y \bigl(z_1, x'_1,e^Y_1
\bigr) \vspace*{2pt}
\cr
\vdots& & & & \vdots\vspace*{2pt}
\cr
z_{m} &
x_{m} & x_m' & e^Y_m
& f_Y \bigl(z_m, x'_m,e^Y_m
\bigr) } \rright).
\end{equation}
A simple method to simulate the i.i.d. copy is to apply some random
permutation $\pi\in S_m$ to $x_1,\ldots,x_n$ and
obtain $x_{\pi(1)},\ldots,x_{\pi(n)}$ (see \cite
{suppcausalstrength}, S.1).
Deleting several arrows with source node $X$ requires several identical
copies $X',X'',\ldots$ of $X$, each generated by
a different permutation.

We then
throw away the two noise columns, that is, the original noise $E^Y$ and
the additional noise $X'$:
%
%
\begin{equation}
\label{eeeedint} \lleft(\matrix{ Z & X & Y \vspace*{2pt}
\cr
\hline
z_1 & x_1 & f_Y \bigl(z_1,
x'_1,e^Y_1\bigr)
\vspace*{2pt}
\cr
\vdots& & \vdots\vspace*{2pt}
\cr
z_m &
x_m & f_Y \bigl(z_m, x'_m,e^Y_m
\bigr) } \rright).
\end{equation}

To see that this triple is indeed sampled from the desired distribution
$P_S(X,Y, Z)$, we recall that
the original structural equation simulates the conditional $P(Y|X,Z)$.
After inserting $X'$ we obtain
the new conditional $\sum_{x'} P(Y|x',Z)\times P(x')$. Multiplying it with
$P(X,Z)$ yields $P_S(X,Y,Z)$, by definition.
Using the above samples from $P_S(X,Y,Z)$ and samples from $P(X,Y,Z)$
we can estimate
\[
\ci_{X\rightarrow Y}=D\bigl(P(X,Y,Z)\|P_S(X,Y,Z)\bigr)
\]
using some known schemes \cite{Perez-Cruz} for estimating relative
entropies from empirical data.
It is important that the samples from the two distributions are
disjoint, meaning that we need to split
the original sample into two halves, one for $P$ and one for~$P_S$.
\end{Ex}

The generation of $P_S$ for a set $S$ of arrows works similarly: every
input of a structural equation that corresponds to an
arrow to be removed is fed with an independent copy of the respective variable.
Although it is conceptually simple to estimate causal strength by
generating the entire joint distribution $P_S$, Theorem~\ref{thmsets}(a)
will show how to break the problem into parts that make estimation of
relative entropies from finite data more feasible.

We now revisit mediation analysis \cite{Pearl00,Avin,Robins}, which
is also based on structural equations, and mention
an interesting relation to our work. Although we have pointed out that
intervening by ``cutting edges'' is complementary to the
intervention on nodes considered there, distributions like $P_S$ can
also occur in an implicit way. To explore the indirect effect
$X\rightarrow Z \rightarrow Y$ in Figure~\ref{DAGs}(b), one can study the
effect of $X$
on $Y$ in the reduced DAG $X\rightarrow Z\rightarrow Y$ under the
distribution $P_{X\rightarrow Y}$ or under
the distribution obtained by setting the copy $X'$ to some fixed value $x'$.
Remarkably, cutting $X\rightarrow Y$ is then used to study the strength
of the \textit{other path} while we use it to
study the strength\footnote{We are grateful to an anonymous referee for
this observation.} of $X\rightarrow Y$.

\subsection{Properties of causal strength}
\label{secproperties}

This subsection shows that our definition of causal strength satisfies
postulates P0--P4. We observe at the same time
some other useful properties.
We start with a property that is used to show
P0.

\textit{Causal strength majorizes observed dependence.}
Recalling that $P(X_1,\ldots,\break X_n)$ factorizes into $\prod_j P(X_j|PA_j)$
with respect to the true causal DAG $G$, one may ask how much error one
would cause
if one was not aware of all causal influences and erroneously assumed
that the true DAG would be the one where some set $S$ of arrows is missing.
The conditionals with respect to the reduced set of parents define a
different joint distribution.

%
\begin{Def}[(Distribution after ignoring arrows)]
\label{dignore}
Given distribution $P$ Markovian with respect to $G$ and set of arrows
$S$, let the partially observed distribution (where interactions across
$S$ are hidden) for node $X_j$ be
\[
\tilde{P}_S\bigl(x_j|pa_j^{\bar{S}}
\bigr) = \sum_{pa_j^S} P\bigl(x_j|pa_j^S,pa_j^{\bar
{S}}
\bigr) P\bigl(pa_j^S|pa_j^{\bar{S}}
\bigr).
\]
Let the partially observed distribution for all the nodes be the product
%
%
\begin{equation}
\label{tildeP} \tilde{P}_S(x_1,\ldots,x_n) =
\prod_j \tilde{P}_S
\bigl(x_j|pa_j^{\bar{S}}\bigr).
\end{equation}
\end{Def}

%
\begin{Rem}
Intuitively, the observed influence of a set of arrows should be
quantified by comparing the data available to an observer who can see
the entire DAG with the data available to an observer who sees all the
nodes of the graph, but only some of the arrows. Definition~\ref
{dignore} formalizes ``seeing only some of the arrows.''
\end{Rem}

Building on Remark~\ref{observeintervene}, the definition of the \emph
{observed dependence} of a set of arrows takes the same general form as
for causal influence. However, instead of inserting noise on the
arrows, we instead simply prevent ourselves from seeing them.

%
\begin{Def}[(Observed influence)]
Given a distribution $P$ that is Markovian with respect to $G$ and set of
arrows\vadjust{\goodbreak}
$S$, let the observed influence of the arrows in $S$ be
\[
\label{empe} \oi_S(P):= D(P\|\tilde{P}_S),
\]
with $\tilde{P}_S$ defined in (\ref{tildeP}).
\end{Def}

The following result, proved in Section~\ref{proofmaj},
is crucial to proving P0.

%
\begin{Thm}[(Causal influence majorizes observed dependence)]
\label{thmgeneralmaj}
Causal influence decomposes into observed influence plus a \emph
{nonnegative} term quantifying the divergence between the partially
observed and interventional distributions
%
%
\begin{equation}
\label{intobs} \ci_S(P) = \oi_S(P) + \sum
_{j=1}^n P\bigl(pa_j^{\bar{S}}
\bigr)\cdot D \bigl(\tilde{P}_S\bigl(X_j|pa_j^{\bar{S}}\bigr)
\| P_S\bigl(X_j|pa_j^{\bar{S}}
\bigr)\bigr).
\end{equation}
\end{Thm}

The theorem shows that ``snapping upstream dependencies'' by using
purely local data that is, by marginalizing using the distribution of
the source node $P(X_i)$ rather than the conditional $P(X_i|PA_i)$---is
essential to quantifying causal influence.


\emph{Proof of postulates for causal strength.}

P0:
Let $G_S$ be the DAG obtained by removing the arrows in $S$ from $G$.
Let $PA^{\bar{S}}_j$ be the parents of $X_j$ in $G_S$, that is, those
that are not in $S$ and introduce the set of nodes
$Z_j$ such that $PA_j =PA_j^{\bar{S}} \cup Z_j$.
By Theorem~\ref{thmgeneralmaj}, $\ci_S=0$ implies $\oi_S=0$, that is,
$\tilde{P}_S=P$, which implies
%
%
\begin{equation}
\label{Z} P(X_j|pa_j)=P\bigl(X_j|pa_j^{\bar{S}}
\bigr)\qquad\forall pa_j^{\bar{S}} \mbox{ with } P
\bigl(pa_j^{\bar{S}}\bigr)\neq0,
\end{equation}
that is,
$X_j
\independent Z_j |PA_j^{\bar{S}}$.

We use again the Ordered Markov condition
%
%
\begin{equation}
\label{altMC} X_j \independent PR_j |
PA_j \qquad\forall j,
\end{equation}
where $PR_j$ denote the predecessors of $X_j$ with respect to some
ordering of nodes that is consistent with $G$. By the contraction rule
\cite{Pearl00},
(\ref{Z}) and (\ref{altMC}) yields
\[
X_j \independent PR_j \cup Z_j |
PA_j^{\bar{S}},
\]
and hence
\[
X_j \independent PR_j |PA_j^{\bar{S}}
,
\]
which is the Ordered Markov condition for $G_S$ if we use the same
ordering of nodes for $G_S$.

P1:
One easily checks $\ci_{X\rightarrow Y}=I(X;Y)$ for the 2-node DAG
$X\rightarrow Y$, because
$
P_{X\rightarrow Y} (x,y)=P(x)P(y),
$
and thus
\[
D(P\|P_{X\rightarrow Y})=D\bigl(P(X,Y)\|P(X)P(Y)\bigr)=I(X;Y).
\]

P2:
Follows from the following lemma.

%
\begin{Lem}[(Causal strength as local relative entropy)]\label{lemdown}
Causal strength $\ci_{k\rightarrow l}$ can be written as the following
relative entropy distance or \textit{conditional} relative entropy distance:
\begin{eqnarray*}
\ci_{k\rightarrow l}
&=& \sum_{pa_l} D \bigl[P(X_l|pa_l)
\| P_S(X_l|pa_l) \bigr]
P(pa_l)\\
&=& D \bigl[ P(X_l|PA_l) \|
P_S(X_l|PA_l) \bigr].
\end{eqnarray*}
\end{Lem}

Note that $P_S(X_l|pa_l)$ actually depends on the reduced set of
parents $PA_l \setminus X_k$ only, but it is more convenient for the
notation and the proof to
keep the formal dependence on all $PA_l$.

\begin{pf*}{Proof of Lemma \ref{lemdown}}
\begin{eqnarray*}
D(P\|P_S) &=&\sum_{x_1\cdots x_n}
P(x_1\cdots x_n) \log\frac
{P(x_1\cdots x_n)}{P_S(x_1 \cdots x_n)} \\
&=&\sum
_{x_1\cdots x_n} P(x_1\cdots x_n) \log\prod
_{j=1}^n \frac
{P(x_j|pa_j)}{P_S(x_j|pa_j)}
\\
&=&\sum_{j=1}^n \sum
_{x_j,pa_j} P(x_j,pa_j) \log
\frac
{P(x_j|pa_j)}{P_S(x_j|pa_j)} \\
&=&\sum_{j=1}^n D \bigl[
P(X_j|PA_j) \| P_S(X_j|PA_j)
\bigr].
\end{eqnarray*}
%
For all $j\neq l$ we have $D [ P(X_j|PA_j)\| P_S(X_j|PA_j) ]
=0$, because
$P(X_l|PA_l)$ is the only conditional that is modified by the deletion.
\end{pf*}

P3:
Apart from demonstrating the postulated inequality,
the following result shows that we have the equality $\ci
_{X\rightarrow
Y}=I(X;Y |PA_Y^X)$
for independent causes.
To keep notation simple, we have restricted our attention to the case
where $Y$ has only two causes $X$ and $Z$,
but $Z$ can also be interpreted as representing all parents of $Y$
other than $X$.

%
\begin{Thm}[(Decomposition of causal strength)]\label{thmdeco}
For the DAGs in Figure~\ref{DAGs}, we have
%
%
\begin{equation}
\label{deco} \ci_{X\rightarrow Y} =I(X;Y |Z) + D \bigl[P(Y|Z)\|
P_{X\rightarrow
Y}(Y|Z)
\bigr].
\end{equation}
If $X$ and $Z$ are independent, the second term vanishes.
\end{Thm}

\begin{pf}
Equation~(\ref{deco}) follows from Theorem~\ref{thmgeneralmaj}: First,
we observe\break $\oi_S(P)=I(X;Y |Z)$ because both measure
the relative entropy distance\vadjust{\goodbreak} between $P(X,Y,Z)$ and
$\tilde{P}_S(X,Y,Z)=P(X,Z)P(Y|Z)$. Second, we have
\begin{eqnarray*}
P_{S}(X,Y,Z)&=& P(X,Z) P_{X \rightarrow Y}(Y|Z).
\end{eqnarray*}
The second summand in (\ref{intobs}) reduces to
\begin{eqnarray*}
&&\sum_z P(z) D\bigl[\tilde{P}_S(Y|z)\|P_S(Y|z)
\bigr]\\
&&\qquad=\sum_z P(z) D\bigl[{P}_S(Y|z)
\|P(Y|z)\bigr] \\
&&\qquad=D\bigl[P(Y|Z)\|P_S (Y|Z)\bigr].
\end{eqnarray*}
%
To see that the second term in equation (\ref{deco}) vanishes for independent $X,Z$, we observe
$P_{X\rightarrow Y}(Y|Z)=P(Y|Z)$ because
\[
P_{X\rightarrow Y}(y|z)= \sum_x P(y|x,z)P(x) =\sum
_x P(y|x,z)P(x|z)=P(y|z). 
\]
\upqed
\end{pf}

Theorem~\ref{thmdeco} states that conditional mutual information
underestimates causal strength.
Assume, for instance, that $X$ and $Z$ are almost always equal because
$Z$ has such a strong influence on $X$ that it is an almost perfect
copy of it.
Then $I(X;Y |Z)\approx0$ because knowing $Z$ leaves almost no
uncertainty about $X$.
In other words, strong dependencies between the causes $X$ and $Z$
makes the influence of cause $X$ almost invisible
when looking at the conditional mutual information $I(X;Y |Z)$ only.
The second term in (\ref{deco}) corrects for the underestimation.
When $X$ depends deterministically on $Z$, it is even the only
remaining term (here, we have again assumed that the conditional distributions
are defined for events that do not occur in observational data).

To provide a further interpretation of
Theorem~\ref{thmdeco}, we recall that $I(X;Y |Z)$ can be seen as the
impact of ignoring the edge $X\rightarrow Y$; see
Remark~\ref{observeintervene}. Then the impact of cutting
$X\rightarrow Y$ is given by the impact of ignoring this link
plus the impact that cutting has on the
conditional $P(Y|Z)$.

P4: This postulate is part (d) of the following
collection of results that relates strength of sets to its subsets.


%
\begin{Thm}[(Relation between strength of sets and subsets)]\label
{tproperties}
\label{thmsets}
The causal influence given in Definition~\ref{dciset} has the
following properties:

\begin{longlist}[(a)]
\item[(a)] \emph{Additivity regarding targets}.
Given set of arrows $S$, let $S_i=\{s\in S|\operatorname{trg}(s)=X_i\}$, then
\[
\ci_S(P)=\sum_i \ci_{S_i}(P).
\]

\item[(b)] \emph{Locality}.
Every $\ci_{S_i}$ only depends on the conditional $P(X_i|PA_i)$ and the
joint distribution of all parents $P(PA_i)$.

\item[(c)] \emph{Monotonicity}.
Given sets of arrows $S_1\subset S_2$ targeting single node $Z$, such
that the source nodes in $S_1$ are jointly independent
and independent of the other parents of $Z$. Then we have
\[
\ci_{S_1}(P)\leq\ci_{S_2}(P).
\]

\item[(d)] \emph{Heredity property}.
Given sets of arrows $S\subset T$, we have
\[
\ci_{T}(P) = 0\quad\implies\quad\ci_{S}(P)=0.
\]
\end{longlist}
\end{Thm}

The proof is presented in Appendix~\ref{subsecproof}.
The intuitive meaning of these properties is as follows.
Part (a) says that causal influence is additive if the arrows have
different targets.
Otherwise, we can still decompose the set $S$ into equivalence classes
of arrows having the same target and obtain additivity regarding the
decomposition.
This can be helpful for practical applications because
estimating each $D[P(PA_i,X_i)\|P_{S_i}(PA_i,X_i)]$ from empirical data
requires less data then estimating
the distance $D(P\|P_S)$ for the entire high dimensional distributions.

We will show in Section~\ref{subsecparad} that general additivity fails.
Part (b) is an analog of P2 for multiple arrows.
According to (c), the strength of a subset of arrows cannot be smaller
than the strength of its superset, provided that there are no
dependencies among the parent nodes.
Finally, part (d) is exactly our postulate P4.

Parts (c) and (d) suggest that monotonicity may generalize to the case
of dependent parents: $S\subset T\implies\ci_S(P)\leq\ci_T(P)$.
However, the following counterexample due to Bastian Steudel shows this
is not the case.

%
\begin{Ex}[(XOR---counterexample to monotonicity when parents are
dependent)]\label{egxor}
Consider the DAG(a) in Figure~\ref{DAGs} and let the relation between
$X,Y,Z$ be given by the structural equations
%
%
\begin{eqnarray}
X &=& Z, \label{cop}
\\
Y& = & X \oplus Z \label{xor}.
\end{eqnarray}
Let $P(Z=0)=a$ and $P(Z=1)=1-a$. Letting $S=\{Z\rightarrow X\}$ and
$T=\{Z\rightarrow X, X\rightarrow Y\}$ we find that
\begin{eqnarray*}
\ci_S(P) & =& -a\log(a)-(1-a)\log(1-a)\quad\mbox{and}
\\
\ci_T(P) & =& -\log\bigl( a^2 + (1-a)^2
\bigr).
\end{eqnarray*}
For $a\notin\{\frac{1}{2},0,1\}$, strict concavity of the logarithm
implies $\ci_T(P)<\ci_S(P)$.
\end{Ex}

\subsection{Examples and paradoxes}
\label{subsecparad}

\emph{Failure of subadditivity}:
The strength of a set of arrows is not bounded from above by the sum of
strength of the single arrows.
It can even happen that removing one arrow from a set has no impact on
the joint distribution while removing
all of them has significant impact, which occurs in communication
scenarios that use redundancy.

\begin{figure}

\includegraphics{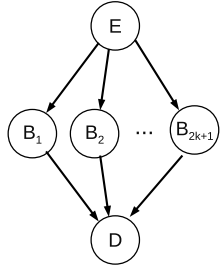}

\caption{Causal structure of an error-correcting scheme: the encoder
generates $2k+1$ bits from a single one.
The decoder decodes the $2k+1$ bit words into a single bit
again.}\label
{figcode}
\end{figure}

%
\begin{Ex}[(Error correcting code)]\label{egerrcode}
Let $E$ and $D$ be binary variables that we call ``encoder'' and
``decoder'' (see Figure~\ref{figcode}) communicating over a channel
that consists of
the bits $B_1,\ldots,B_{2k+1}$. Using the simple repetition code, all
$B_j$ are just copies of $E$. Then $D$ is set to the logical value that
is attained by the
majority of $B_j$.
This way, $k$ errors can be corrected, that is, removing $k$ or less of
the links $B_j\rightarrow D$ has no effect
on the joint distribution, that is,
$P_S=P$ for
$S:=(B_1\rightarrow D,B_2\rightarrow D,\ldots,B_k\rightarrow D)$, hence
$\ci_S(P)=0$. In words: removing $k$ or less arrows is without impact,
but removing all of them
is, of course.
After all, the arrows jointly generate the dependence
$I(E;D)=I((E,B_1,\ldots,B_k);D)=1$, provided that $E$ is uniformly distributed.
\end{Ex}

%

Clearly, the outputs of $E$ causally influence the behavior of $D$.
We therefore need to consider interventions that destroy many arrows
at once if we want to capture the fact that
their joint influence is nonzero.

Thus, causal influence of arrows is \emph{not} subadditive:
the strength of each arrow $B_j\rightarrow D$ is zero, but the strength
of the set of all $B_j\rightarrow D$ is $1$ bit.

\emph{Failure of superadditivity}:
The following example reveals an opposing phenomenon, where the causal
strength of a set is smaller than the sum of the single arrows.

%
\begin{Ex}[(XOR with uniform input)]\label{egxoru}
Consider the structural equations~(\ref{cop}) and (\ref{xor}) with
uniformly distributed $Z$.
The causal influence of each arrow targeting the XOR-gate individually
is the same as the causal influence of both arrows taken together:
\[
\ci_{X\rightarrow Y}(P)=\ci_{Z\rightarrow Y}(P)=\ci_{\{X\rightarrow
Y,Z\rightarrow Y\}}(P)=1 \mbox{ bit.}
\]
\end{Ex}

\emph{Strong influence without dependence/failure of converse of P0}:
Revisiting Example~\ref{egxoru} is also instructive because it
demonstrates an extreme case of confounding where $I(X;Y |Z)$ vanishes
but causal influence is strong.
Removing $X\rightarrow Y$ yields
\[
P_{X\rightarrow Y}(x,y,z)=P(x,z)P(y),
\]
where $P(z)=P(y)=1/2$ and $P(x|z)=\delta_{x,z}$.
It is easy to see that
\[
D(P \| P_{X\rightarrow Y})=1,
\]
because $P$ is a uniform distribution over $2$ possible triples
$(x,y,z)$, whereas $P_{X\rightarrow Y}$ is a uniform distribution over
a superset of $4$ triples.

The impact of cutting the edge $X\rightarrow Y$ is remarkable: both
distributions, the observed one $P$ as well as the post-cutting
distribution $P_S$, factorize
$P_S(X,Y,Z)=P_S(X,Z)P_S(Y)$ and $P(X,Y,Z)=P(X,Z)P(Y)$. Cutting the edge
keeps this product structure and changes the joint distributions by
only changing the marginal distribution of $Y$ from $P(Y)$ to $P_S(Y)$.

Note that $P$ satisfies the Markov condition with respect to
$G_{X\rightarrow Y}$ (i.e., the DAG obtained from the original one
by dropping $X\rightarrow Y$)
because $Y$ is a constant. Since $\ci_{X\rightarrow Y} \neq0$, this
shows that the converse of P0 does not hold.

\begin{figure}

\includegraphics{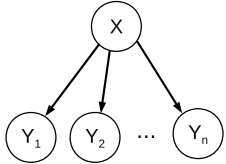}

\caption{Broadcasting one bit from one node to multiple nodes.} \label
{figbroad}
\end{figure}

\emph{Strong effect of little information}:
The following example considers multiple arrows and shows that their
joint strength may even be strong when they carry the same small amount
of information.
%
%
\begin{Ex}[(Broadcasting)]
\label{egbroadcasting}
Consider a single source $X$ with many targets $Y_1,\ldots, Y_n$ such
that each $Y_i$ copies $X$, see Figure~\ref{figbroad}. Assume
$P(X=0)=P(X=1)=\frac{1}{2}$.
If $S$ is the set of all arrows $X\rightarrow Y_j$ then $\ci_S=n$.
Thus, the single node $X$ exerts $n$ bits of causal influence on its dependents.
\end{Ex}
%
%

\section{Causal influence between two time series}

\subsection{Definition}

Since causal analysis of time series is of high practical importance,
we devote a section to this case.
For some fixed $t$, we introduce the short notation $X\rightarrow Y_t$ for\vadjust{\goodbreak}
the set of all arrows that point to $Y_t$ from some $X_s$ with $s<t$.
Then
\[
\ci_{X \rightarrow Y_t}
\]
measures the impact of deleting all these arrows. We propose
to replace transfer entropy with this measure since it does not suffer
from the drawbacks described in Section~\ref{subsecte}.

Section~\ref{subsecse} describes how to estimate causal strength
from finite data for one arrow and briefly mentions how this
generalizes to set of arrows.
To keep this section self-consistent, we briefly rephrase the
description for the case of time series.

Suppose we have learned the structural equation model
%
%
\begin{equation}
\label{structime} Y_t=f_t(X_{t-1},X_{t-2},
\ldots,X_{t-p},E_t),
\end{equation}
from observed data $(x_t,y_t)_{t\leq0}$,
where the noise variables $E_t$ are jointly independent and independent
of $X_t,X_{t-1},\ldots,Y_{t-1},Y_{t-2},\ldots.$
Assume, moreover, that we have inferred the corresponding values
$(e_t)_{t\leq0}$ of the noise.
If we have multiple copies of the time series, we can apply the method
described in Section~\ref{subsecse} in a straightforward way:
Due to the locality property stated in Theorem~\ref{thmsets}(b), we
only consider the variables
$X_{t-p},\ldots,X_{t-1},Y_t$ and feed (\ref{structime}) with i.i.d.
copies of $X_{t-p},\ldots,X_{t-1}$ by applying
random permutations to the observations, which then yields samples from
the modified distribution $P_S(X_{t-p},\ldots,X_{t-1},Y_t)$.

If we have only one observation for each time instance, we have to
assume stationarity (with constant function
$f_t=f$) and ergodicity
and generate an artificial statistical sample by looking at
sufficiently distant windows.

\subsection{Comparison of causal influence with transfer entropy}
\label{subsecTE}

We first recall the example given by \cite{AyInfoFlow} showing a
problem with transfer entropy (Section~\ref{subsecte}).
Assume that the variables $X_t,Y_t$ in Figure~\ref{figTen}, right, are
binary and the transition from $X_{t-1}$ to $Y_t$ is a perfect copy and likewise
the transition from $Y_{t-1}$ to $X_t$.
Assume, moreover, that the two causal chains have been initialized such that, with
probability $1/2$, all variables are $1$ and with probability $1/2$ all
are zero.
Then the set $X\rightarrow Y_t$ is the singleton $S:=\{
X_{t-1}\rightarrow Y_t\}$.
Using Lemma~\ref{lemdown}, we have
\[
\ci_{X_{t-1} \rightarrow Y_t}=D \bigl[P(Y_t,X_{t-1})\|
P_S(Y_t,X_{t-1}) \bigr].
\]
Since $Y_t$ is a perfect copy of $X_{t-1}$, we have
\[
P(y_t,x_{t-1})=\cases{ %
1/2,
& \quad$\mbox{for } x_{t-1}=y_t$,
\vspace*{2pt}\cr
0, & \quad$\mbox{otherwise}$}
\]
into
\[
P_S(y_t,x_{t-1})=1/4 \qquad\mbox{for }
(y_t,x_{t-1}) \in\{0,1\}^2.
\]
One easily checks $D(P\|P_S)=1$.\eject

Note that the example is somewhat unfair, since it is \emph
{impossible} to distinguish the structural equations from a model
without interaction between $X$ and $Y$,
where $X_{t+1}$ is
obtained from $X_t$ by inversion and similarly for $Y$, no matter how
many observations are performed. Thus, from observing the system it is
impossible to tell whether or not $X$ exerts an influence on $Y$.
However, the following modification shows that transfer entropy still
goes quantitatively wrong if
small errors are introduced.

%
\begin{Ex}[(Perturbed transfer entropy counterexample)]\label
{experturbed}
Perturb Ay and Polani's example by having $Y_{t}$ copy $X_{t-1}$
correctly with probability $p=1-\varepsilon$. Set node $Y_t$'s transitions
as Markov matrix
\[
\lleft(\matrix{ & \vline& x_{t-1}=0 & x_{t-1}=1\vspace*{2pt}
\cr
\hline y_{t}=0 & \vline& 1-\varepsilon& \varepsilon\vspace*{2pt}
\cr
y_{t}=1 & \vline& \varepsilon& 1-\varepsilon} \rright),
\]
and similarly for the transition from $Y_{t-1}$ to $X_t$.

The transfer entropy from $X$ to $Y$ at time $t$ is
\begin{eqnarray*}
\mathrm{TE} &=&I(X_{(-\infty,t-1]};Y_t |Y_{(-\infty,t-1]})
=I(X_{t-1};Y_t|Y_{t-2})
\\
&=& H(Y_t|Y_{t-2}) - H(Y_t|Y_{t-2},X_{t-1})=H(Y_t|Y_{t-2})
- H(Y_t|X_{t-1}),
\end{eqnarray*}
where $H(\cdot|\cdot)$ denotes the conditional Shannon entropy.
The equalities can be derived from d-separation in the causal DAG
Figure~\ref{figTen}, right \cite{Pearl00}.
For instance, conditioning on $Y_{t-2}$, renders the pair $(Y_t,X_{t-1})$
independent of all the remaining past of $X$ and $Y$.
We find
\begin{eqnarray*}
-H(Y_t|X_{t-1}) & =& \varepsilon\log\varepsilon+(1-
\varepsilon)\log(1-\varepsilon),
\\
H(Y_t|Y_{t-2})&=&2\varepsilon(1-\varepsilon)\log
\frac{1}{2\varepsilon(1-\varepsilon
)}+\bigl(1-2\varepsilon+2\varepsilon^2\bigr)\log
\frac{1}{1-2\varepsilon+2\varepsilon^2}.
\end{eqnarray*}
Hence,
\begin{eqnarray*}
\mathrm{TE} &=& \bigl(1-2\varepsilon+2\varepsilon^2\bigr)\log\frac
{1}{1-2\varepsilon+2\varepsilon^2}
+ 2\varepsilon(1-\varepsilon)\log\frac{1}{2\varepsilon
(1-\varepsilon)}\\
&&{} + \varepsilon\log
\varepsilon+(1-\varepsilon)\log(1-\varepsilon),
\end{eqnarray*}
which tends to zero as $\varepsilon\rightarrow0$.

Causal influence, on the other hand, is given by the mutual information
$I(Y_t;X_{t-1})$ because all
edges other than $X_{t-1}\rightarrow Y_t$ are irrelevant (see
Postulate~P2). Thus,
\begin{eqnarray*}
\ci_{X\rightarrow Y_t}=H(Y_t)-H(Y_t|X_{t-1}) &=
1+ (1-\varepsilon)\log(1-\varepsilon)+\varepsilon\log\varepsilon,
\end{eqnarray*}
which tends to $1$ for $\varepsilon\to0$.
Hence, causal influence detects the causal interactions between $X$ and
$Y$ based on \emph{empirical data}, whereas transfer entropy does not.
Thanks to the perturbation, the joint distribution tells us the kind of
causal relations by which it is generated. For large enough samples,
the strong discrepancy between
transfer entropy and our causal strength thus becomes apparent.
\end{Ex}

\section{Causal strength for linear structural equations}
\label{seclin}

For linear structural equations, we can provide a more explicit
expression of causal strength under the assumption of multivariate Gaussianity.
Let $n$ random variables $X_1,\ldots,X_n$ be ordered such that there are
only arrows from $X_i$ to $X_j$ for $i<j$.
Then we have structural equations
\[
X_{j} =\sum_{i<j} A_{ij}
X_i +E_j,
\]
where all $E_j$ are jointly independent noise variables.
In vector and matrix notation we have
%
%
\begin{equation}
\label{strma} X=AX+E\qquad\mbox{that is, } X=(I-A)^{-1} E,
\end{equation}
where $A$ is lower triangular with zeros in the diagonal.

To compute the strength of $S\subset\{1,\ldots,n\}$,
we assume for reasons of convenience that all variables have zero mean.
Then $D(P\|P_S)$
can be computed from the covariance matrices alone.

The covariance matrix of $X$ reads
\[
\Sigma= (I-A)^{-1} \Sigma_E (I-A)^{-T},
\]
where $\Sigma_E$ denotes the covariance matrix of the noise (which is
diagonal by assumption) and
$(\cdot)^{-T}$ the transpose of the inverse of a matrix.

To compute the covariance matrix $\Sigma^S$ of $P_S$,
we first split $A$ into $A_S+A_{\bar{S}}$, where $A_S$ contains only
those entries that correspond to the edges in the set $S$
and $A_{\bar{S}}$ only those corresponding to the complement of $S$.
Using this notation, the modified structural equations read
%
%
\begin{equation}
\label{modstr} X= A_{\bar{S}} X +E + A_{S} X',
\end{equation}
where $X'=(X_1',\ldots,X_n')^T$ and each $X'_j$ has the same
distribution as $X_j$ and satisfies joint independence of all
$X'_1,\ldots,X_n', E_1,\ldots,E_n$. It is convenient to define the
modified noise
\[
E':=E+A_S X',
\]
with covariance matrix
%
%
\begin{equation}
\label{primE} \Sigma_{E'}= \Sigma_E + A_S
\Sigma_X^D A_S^T,
\end{equation}
where $\Sigma_X^D$ contains only the diagonal entries of $\Sigma_X$
(recall that all $X_j'$ are independent).
The modified variables $X^S$ are now given by the equation
\[
X^S=A_{\bar{S}} X +E',
\]
which formally looks like (\ref{strma}), although
the components of $E'$ are dependent while the $E_j$ in (\ref{strma})
are independent.\vadjust{\goodbreak}
Thus, we obtain
the modified covariance matrix of $X$ by
\[
\Sigma_S = (I-A_{\bar{S}})^{-1} \Sigma_{E'}
(I-A_{\bar{S}})^{-T}.
\]
The causal strength now reads
\begin{eqnarray*}
\ci_S= D(P\|P_S)&=& \frac{1}{2} \biggl( \operatorname{tr}
\bigl[\Sigma_S^{-1} \Sigma\bigr] -\log\frac{\det\Sigma}{\det
\Sigma_S}-n
\biggr)
\\
&=& \frac{1}{2} \biggl( \operatorname{tr} \bigl[ (I-A_{\bar{S}})
\Sigma_{E'}^{-1} (I-A_{\bar{S}}) (I-A)^{-1}
\Sigma_E (I-A)^{-1} \bigr]
\\
&&\hspace*{54pt}{}-\log\frac{\det(I-A)^{-1}\Sigma_E(I-A)^{-1}}{\det(I-A_{\bar
{S}})^{-1}\Sigma_{E'}(I-A_{\bar{S}})^{-1}}-n \biggr),
\end{eqnarray*}
with $\Sigma_{E'}$ given by (\ref{primE}).



%
\begin{Ex}[(Linear structural equations with independent
parents)]\label
{exlinear}
It is instructive to look at the following simple case:
\[
X_n:= \sum_j \alpha_{nj}
X_j +E_n\qquad\mbox{with } E_n,X_1,
\ldots,X_{n-1} \mbox{ jointly independent.}
\]
For the set
$S:=\{X_1\rightarrow X_n,\ldots,X_k\rightarrow X_n\}$
with $k\leq n$ some calculations show
\[
\ci_S = \frac{1}{2}\log\frac{\operatorname{Var} (X_n) -\sum_{j=k+1}^{n-1}
\alpha_{nj}^2 \operatorname{Var} (X_j) }{\operatorname{Var} (X_n)
-\sum_{j=1}^{n-1}
\alpha_{nj}^2 \operatorname{Var} (X_j)}.
\]
For the single arrow $X_1\rightarrow X_n$, we thus obtain
\[
\ci_{X_1 \rightarrow X_n} = \frac{1}{2} \log\frac{\operatorname
{Var} (X_n) -\sum_{j=2}^{n-1}
\alpha_{nj}^2 \operatorname{Var} (X_j) }{\operatorname{Var} (X_n)
-\sum_{j=1}^{n-1}
\alpha_{nj}^2 \operatorname{Var} (X_j)}.
\]
If $X_1$ is the only parent, that is, $n=2$, we have
\[
\ci_{X_1 \rightarrow X_2} = \frac{1}{2}\log\frac{\operatorname
{Var} (X_2)}{\operatorname{Var} (X_2) -
\alpha_{21}^2 \operatorname{Var} (X_1)} = -
\frac{1}{2} \log(1-r_{21}),
\]
with $r_{21}$ as in equation~(\ref{rdef}) introduced in the context
of ANOVA. Note that the relation between our measure and $r_{n1}$ is
less simple for $n>2$ because
$r_{n1}$ would then still measure the fraction of the variance of $X_n$
explained by $X_1$, while $\ci_{X_1\rightarrow X_n}$
is related to the fraction of the \emph{conditional} variance of $X_n$,
given its other parents, explained by $X_1$.
This is because our causal strength reduces to a \emph{conditional}
mutual information for independent parents; see the last sentence of
Theorem~\ref{thmdeco}.
\end{Ex}

\section{Experiments}

Code for all experiments can be downloaded at\break
\url{http://webdav.tuebingen.mpg.de/causality/}.

\subsection{DAGs without time structure}

We here restrict attention to linear structural equations, but
interesting generalizations are given by additive noise models
\cite{Hoyer,UAIidentifiability,Jonastpami} and post-nonlinear models
\cite{ZhangUAI}.

The first step in estimating the causal strength consists in inferring
the structure matrix $A$ in (\ref{strma}) from
the given matrix $\mathbf{X}$ of observations $x_j^i$ with $j=1,\ldots,n$
and $i=1,\ldots,2k$
(the $j$th row corresponds to the observed values of $X_j$).
We did this step by ridge regression. We decompose $A$ into the sum
$A_S+A_{\bar{S}}$ as in Section~\ref{seclin}.

Then we divide the columns of $\mathbf{X}$ into two parts $\mathbf
{X}_A$ and
$\mathbf{X}_B$ of sample size $k$.
While $\mathbf{X}_A$ is kept as it is, $\mathbf{X}_B$ is used to
generate new
samples according to the modified structural equations:
First, we note that the values of the noise variables corresponding to
the observations $\mathbf{X}_B$ are given by
the residuals
\[
\mathbf{E}_B:=\mathbf{X}_B -A\cdot
\mathbf{X}_B.
\]
Then we generate a matrix $\mathbf{X'_B}$ by applying independent random
permutations to the columns of $\mathbf{X}_B$, which simulates samples of
the random variables $X_j'$ in (\ref{modstr}).
Samples from the modified structural equation are now given by
\[
\mathbf{X}_B^S:=(I-A_{\bar{S}})^{-1}
\cdot\mathbf{X}_B + \mathbf{E}_B + A_S \cdot
\mathbf{X}'_B.
\]

\begin{figure}

\includegraphics{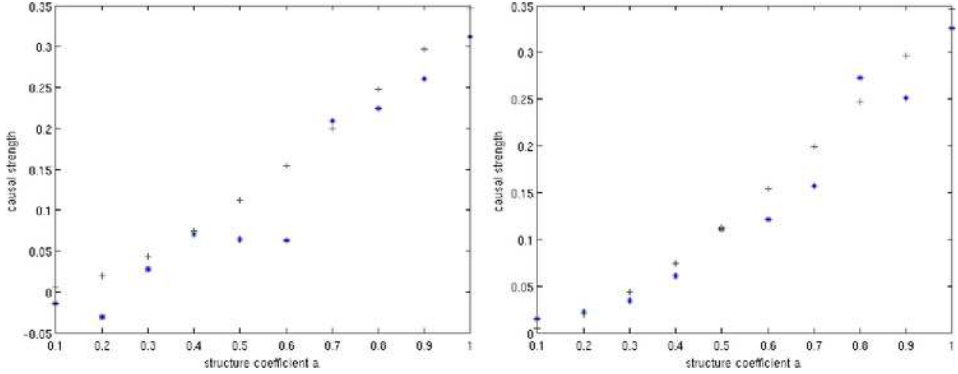}

\caption{Estimated and computed value $\ci_{1
\rightarrow2}$ for $X_1 \rightarrow X_2$, indicated by $*$ and $+$,
respectively.
The underlying linear Gaussian model reads $X_2=a\cdot X_1+E$. Left for
sample size $1000$, which amounts to $500$ samples in each part.
Right: sample size $2000$, which yields more reliable results.}
\label{figbiv}
\end{figure}

To estimate the relative entropy distance between $P$ and $P_S$ (with
samples $\mathbf{X}_A$ and $\mathbf{X}^S_B$),
we use the method described in \cite{Perez-Cruz}:
Let $d_i$ be the euclidean distance from the $i$th column in ${\bf
X}_A$ to the $r$th nearest neighbor among the other columns of ${\bf
X}_A$ and
$d_i^S$ be the distance to the $r$th nearest neighbor among all columns
of $\mathbf{X}_B$, then the estimator reads
\[
\hat{D}(P\|P_S):=\frac{n}{k} \sum
_{i=1}^k \log\frac{d^S_i}{d_i} +\log
\frac{k}{k-1}.
\]
Figure~\ref{figbiv} shows the difference between estimated and
computed causal strength for the simplest
DAG $X_1\rightarrow X_2$ with increasing structure coefficient.
%
For some edges, we obtain significant bias. However,
since the bias depends on the distributions~\cite{Perez-Cruz}, it would
be challenging to correct for it.

%
%

%
%
\begin{figure}[b]

\includegraphics{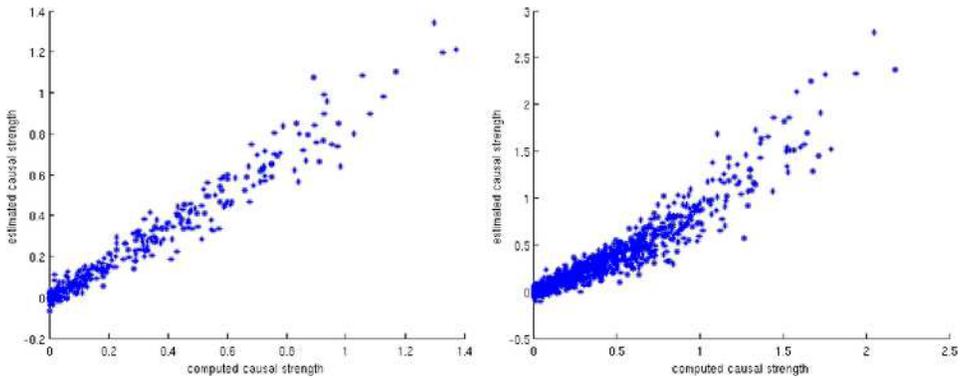}

\caption{Relation between computed and estimated single arrow strengths
for $100$ randomly generated structure matrices and noise variance $1$.
The estimation is based on sample size $1000$. Left: complete DAG on
$3$ nodes.
Right: the same for $5$ nodes.}\label{figscatter}
\end{figure}

To provide a more general impression on the estimation error, we have
considered a complete DAG on $n=3$ and $n=6$ nodes and randomly
generated structure coefficients.
In each of $\ell=1,\ldots,100$ runs, the structure matrix is generated
by independently
drawing each entry from a standard normal distribution. For each of the
${n \choose2}$ arrows $i\rightarrow j$ and each $\ell$ we computed and
estimated
$\ci_{i\rightarrow j}$, which yields the $x$-value and the $y$-value,
respectively, of one of the $100 \cdot{n \choose2}$ points in the
scatter plots in
Figure~\ref{figscatter}.
Remarkably, we do not see a significant degradation for $n=6$ nodes
(right) compared to $n=3$ (left).


\subsection{Time series}

The fact that transfer entropy fails to capture causal strength has
been one of our motivations for
defining a different measure. We revisit the critical example in
Section~\ref{subsecTE},\vadjust{\goodbreak} where the dynamical evolution on two bits
was given by noisy copy operations from
$X_{t-1}$ to $Y_t$ and $Y_{t-1}$ to $X_t$. This way, we obtained causal
strength $1$ bit when the copy operations is getting perfect.
Our software for estimating causal strength only covers the case of
linear structural equations, with the additional assumption of
Gaussianity for the subroutines that \textit{compute} the causal strength
from the covariance matrices for comparison with the estimated value.

A natural linear version of Example~\ref{experturbed} is an
autoregressive (AR-) model of order $1$ given by
\[
\pmatrix{ X_t
\cr
Y_t }= \pmatrix{ 0 & \sqrt{1-\varepsilon^2}
\cr
\sqrt{1-\varepsilon^2} & 0 }
\pmatrix{ X_{t-1}
\cr
Y_{t-1} } + \pmatrix{ E^X_t
\cr
E^Y_t },
\]
where $E^X_t,E^Y_t$ are independent noise terms. We consider the
stationary regime where $X_t$ and $Y_t$ have unit variance and
$E_t$ has variance $\varepsilon^2$.
For $\varepsilon\to0$ the influence from $X_{t-1}$ on $Y_t$, and
similarly from $Y_{t-1}$ to $X_t$ gets deterministic. We thus obtain infinite
causal strength (note that two deterministically coupled random
variables with probability density have infinite mutual information).
It is easy to see that transfer entropy does not diverge, because the
conditional variance of $Y_t$ is $2\varepsilon^2$
if only the past of $Y$ is given and $\varepsilon^2$ if the past of
$X$ is
given in addition. Reducing the variance by the factor $2$ corresponds to
the constant information gain of $\frac{1}{2} \log2$, regardless of
how small $\varepsilon$ is.

\begin{figure}

\includegraphics{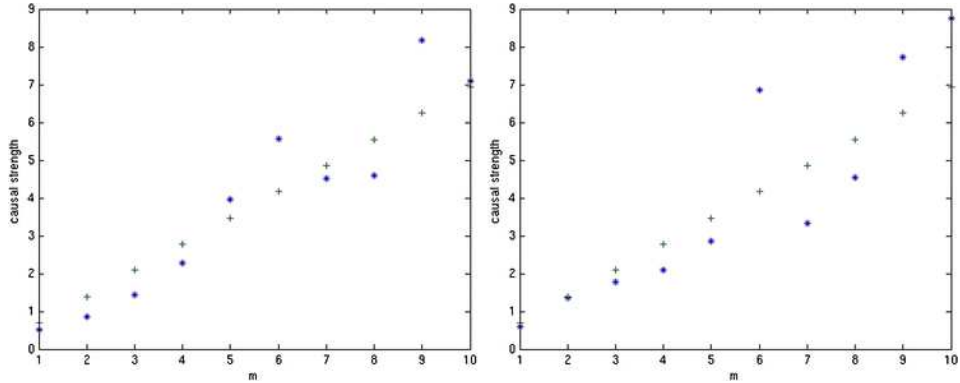}

\caption{Estimated and computed value $\ci_{X \rightarrow Y_t}$ where
$\varepsilon=2^{-m}$
and $m$ runs from $1$ to $10$. Left: for length $T=5000$. Right:
$T=50\mbox{,}000$.}\label{figtscopy}
\end{figure}

Figure~\ref{figtscopy} shows the computed and estimated values of
causal strength for decreasing $\varepsilon$, that is, the
deterministic limit.
Note that, in this limit, the estimated relative entropy can deviate
strongly from the true one because the true one diverges since
$P_S$ lives on a higher dimensional manifold than $P$. This probably
explains the large errors for $m\geq 6$, which correspond to
quite low noise level already.
%

%

\section{Conclusions}

We have defined the strength of an arrow or a set of arrows in a
causal Bayesian network by quantifying the impact of an operation that
we called ``destruction of edges''. We have stated a few postulates that
we consider natural for a measure of causal strength and
shown that they are satisfied by our measure. We do not claim that our
list is complete, nor do we claim that
measures violating our postulates are inappropriate.
How to quantify causal influence may strongly depend on the purpose
of the respective measure.

For a brief discussion of an alternative measure of causal strength and
some of the difficulties that arising when quantifying the total
influence of one set of nodes on another, see the supplementary
material \cite{suppcausalstrength}.

The goal of this paper is to encourage discussions on
how to define causal strength within a framework that is general enough to
include dependencies between variables of arbitrary domains, including
nonlinear interactions, and multi-dimensional and discrete variables
at the same time.

\begin{appendix}\label{app}
\section*{Appendix: Further properties of causal strength and proofs}

\subsection{\texorpdfstring{Proof of Theorem~\protect\ref{thmgeneralmaj}}{Proof of Theorem 2}}

\label{proofmaj}

Expand $\ci_S(P)$ as
%
%
\begin{eqnarray}
\label{KLsum} D (P\|P_S ) & =& \sum_{x_1\cdots x_n}
P(x_1\cdots x_n)\log\frac
{P(x_1\cdots x_n)}{P_S(x_1\cdots x_n)}
\\
& =& \sum_{x_1\cdots x_n} P(x_1\cdots
x_n)\log\frac{P(x_1\cdots
x_n)}{\tilde{P}_S(x_1\cdots x_n)}
\nonumber
\\[-8pt]
\\[-8pt]
\nonumber
&&{}+ \sum_{x_1\cdots x_n}
P(x_1\cdots x_n)\log\frac{\tilde{P}_S(x_1\cdots x_n)}{P_S(x_1\cdots x_n)}.
\end{eqnarray}

Note that the second term can be written as
%
%
\begin{eqnarray}
\label{KLsumS} &&\sum_{x_1\cdots x_n} P(x_1\cdots
x_n)\log\prod_{j=1}^n
\frac{\tilde
{P}_S(x_j|pa_j^{\bar{S}})}{P_S(x_j|pa_j^{\bar{S}})}
\nonumber
\\[-8pt]
\\[-8pt]
\nonumber
&&\qquad= \sum_{j=1}^n \sum
_{x_1\cdots x_n} P(x_1\cdots x_n)\log
\frac{\tilde
{P}_S(x_j|pa_j^{\bar{S}})}{P_S(x_j|pa_j^{\bar{S}})}
\\
&&\qquad= \sum_{j=1}^n \sum
_{x_j,pa_j} P\bigl(x_j,pa^{\bar{S}}_j,
pa^S_j\bigr)\log\frac
{\tilde{P}_S(x_j|pa_j^{\bar{S}})}{P_S(x_j|pa_j^{\bar{S}})}
\\
&&\qquad= \sum_{j=1}^n \sum
_{x_j,pa_j^{\bar{S}}} \tilde{P}\bigl(x_j|pa^{\bar
{S}}_j
\bigr)P\bigl(pa^{\bar{S}}_j\bigr)\log\frac{\tilde
{P}_S(x_j|pa_j^{\bar
{S}})}{P_S(x_j|pa_j^{\bar{S}})}
\\
&&\qquad= \sum_{j=1}^n P\bigl(pa_j^{\bar{S}}
\bigr)\cdot D \bigl[\tilde{P}_S\bigl(X_j|pa_j^{\bar
{S}}\bigr)
\| P_S\bigl(X_j|pa_j^{\bar{S}}\bigr)
\bigr].
\end{eqnarray}





Causal influence is thus observed influence plus a correction term that
quantifies the divergence between the partially observed and
interventional distributions. The correction term is nonnegative since
it is a weighted sum of conditional Kullback--Leibler divergences.

\subsection{Decomposition into conditional relative entropies}

The following result generalizes Lemma~\ref{lemdown} to the case where
$S$ contains more than one edge.\vadjust{\goodbreak}
It shows that the relative entropy expression defining causal strength
decomposes into a sum of conditional relative entropies, each of it
referring to the conditional distribution of one of the target nodes,
given its parents:

%
%
\begin{Lem}[(Causal influence decomposes into a sum of expectations)]%
\label{tlocaladd}
The causal influence of set of arrows $S$ can be rewritten
\[
\ci_S(P) = \sum_{j\in \operatorname{trg}(S)} D
\biggl(P(X_j|PA_j) \Big\Vert\sum
_{pa_j^S} P\bigl(X_j|PA^{\bar{S}}_j,pa^S_j
\bigr)\cdot P_{\prod}\bigl(pa^S_j\bigr) \biggr), \label{elocal}
\]
where $\operatorname{trg}(S)$ denotes the target nodes of arrows in $S$.
\end{Lem}

The result is used in the proof of Theorem~\ref{thmsets} below.


\begin{pf*}{Proof of Theorem \ref{thmsets}}
Using the chain rule for relative entropy \cite{cover}, we get
%
%
\begin{eqnarray}
 D(P\|P_S)&=& \sum_{j=1}^n
D \bigl[P(X_j|PA_j) \Vert P_S(X_j|PA_j)
\bigr]\\
\label{elocal1}&=&
\sum_{j=1}^n \sum
_{pa_j} P(pa_j) D \bigl[P(X_j|pa_j)
\Vert P_S(X_j|pa_j) \bigr]
\nonumber
\\[-8pt]
\\[-8pt]
\nonumber
& =&\sum
_{j\in \operatorname{trg}(S)} D \bigl[P(X_j|PA_j)
\Vert P_S(X_j|PA_j) \bigr],
\end{eqnarray}
where we have used that $P(X_j|PA_j)=P_S(X_j|PA_j)$ for all $j\notin
\operatorname{trg}(S)$. Then the statement follows from the definition of
$P_S(X_j|PA_j)$. Note that a similar statement for $D(P_S\|P)$ (i.e.,
swapping the roles of $P$ and $P_S$) would not hold because then the
weighting factor $P(pa_j)$ in (\ref{elocal1}) needed to be replaced
with the factor
$P_S(pa_j)$, which is sensitive even to deleting edges not targeting $j$.
\end{pf*}

\subsection{\texorpdfstring{Proof of Theorem~\protect\ref{thmsets}}{Proof of Theorem 5}}
\label{subsecproof}

Parts (a) and (b) follow from Lemma~\ref{tlocaladd} since $\ci
_{S_i}(P)$ is the $i$th summand in (\ref{elocal1}), which obviously depends
on $P(X_i|PA_i)$ and $P(PA_i)$ only.

To prove part (c), we will show that the restrictions
of $P,P_{S_1},P_{S_2}$ to the variables $Z,PA_Z$ form a so-called
Pythagorean triple in the sense of \cite{Amari}, that is,
%
%
\begin{eqnarray}
\label{orth} \qquad & &D \bigl[P(Z,PA_Z) \| P_{S_2}(Z,PA_Z)
\bigr]
\nonumber
\\[-8pt]
\\[-8pt]
\nonumber
&&\qquad=
D \bigl[P(Z,PA_Z) \| P_{S_1}(Z,PA_Z) \bigr]
\nonumber
+ D\bigl[P_{S_1}(Z,PA_Z)\| P_{S_2}(Z,PA_Z)
\bigr].
\end{eqnarray}
This is sufficient
because the left-hand side and the first term on the right-hand side of
equation~(\ref{orth}) coincide with $\ci_{S_2}$ and $\ci_{S_1}$,
respectively,
due to
part (b). Note, however, that
\[
D \bigl[P_{S_1}(Z,PA_Z) \| P_{S_2}(Z,PA_Z)
\bigr]\neq D (P_{S_1} \| P_{S_2} )
\]
because
we have such a locality statement only for terms of the form $D(P\|
P_S)$. We therefore
consistently restrict attention to $Z,PA_Z$ and
find
\begin{eqnarray*}
&&D \bigl[P(Z,PA_Z) \| P_{S_2}(Z,PA_Z) \bigr]
\\
&&\qquad= \sum_{z,pa_Z} P(z,pa_Z) \log
\frac{P(z|pa_Z)}{P_{S_2}(z|pa^{\bar{S}_2}_Z)}
\\
&&\qquad= \sum_{z,pa_Z} P(z,pa_Z) \log
\frac{P(z|pa_Z)}{P_{S_1}(z|pa^{\bar
{S}_1}_Z)} + \sum_{z,pa_Z} P(z,pa_Z)
\log\frac{P_{S_1}(z|pa^{\bar
{S}_1}_Z)}{P_{S_2}(z|pa^{\bar{S}_2}_Z)}
\\
&&\qquad= D \bigl[P(Z,PA_Z) \| P_{S_1}(Z,PA_Z)
\bigr] \\
&&\qquad\quad{}+ \sum_{z,pa_Z} P\bigl(z|pa^{\bar{S}_1}_Z,pa^{S_1}_Z
\bigr) P_{\prod
}\bigl(pa^{S_1}_Z\bigr)P
\bigl(pa^{\bar{S}_1}_Z\bigr) \log\frac{P_{S_1}(z|pa^{\bar
{S}_1}_Z)}{P_{S_2}(z|pa^{\bar{S}_2}_Z)},
\end{eqnarray*}
where we have
used that the sources in $S_1$ are
jointly independent and independent of the other parents of $Z$.
By definition of $P_{S_1}$, the second summand reads
\[
\sum_{z,pa^{\bar{S}_1}_Z} P_{S_1} \bigl(z,pa^{\bar{S}_1}_Z
\bigr) \log\frac{P_{S_1}(z|pa^{\bar{S}_1}_Z)}{P_{S_2}(z|pa^{\bar{S}_2}_Z)}
\\
= D \bigl[P_{S_1}(Z,PA_Z)
\| P_{S_2}(Z,PA_Z) \bigr],
\]
which proves (\ref{orth}).

By Lemma~\ref{tlocaladd}, it is only necessary to prove part (d) in
the case where both $S$ and $T$ consist of arrows targeting a single
node. To keep the exposition simple, we consider the particular case of
a DAG containing three nodes $X,Y,Z$ where $S=\{X\rightarrow Z\}$ and
$T=\{X\rightarrow Z, Y\rightarrow Z\}$. The more general case follows
similarly. Observe that $D(P\|P_T)=0$ if and only if
%
%
\begin{equation}
\label{econd} P(Z|x,y)=\sum_{\hat{x},\hat{y}}P(Z|\hat{x},
\hat{y})P(\hat{x})P(\hat{y})
\end{equation}
%
for all $x,y$ such that $P(x,y)>0$. Multiplying both sides with $P(x')$
and summing over all $x'$
yields
\[
\sum_{x'} P\bigl(Z|x',y\bigr)P
\bigl(x'\bigr)=\sum_{\hat{x},\hat{y}}P(Z|\hat{x},
\hat{y})P(\hat{x})P(\hat{y}),
\]
because the right-hand side does not depend on $x$. Using (\ref
{econd}) again, we obtain
\[
\sum_{x'} P\bigl(Z|x',y\bigr)P
\bigl(x'\bigr)=P(Z|x,y)
\]
for all $x,y$ with $P(x,y)\neq0$. Hence $P_S=P$, and thus $D(P\|P_S)=0$.

\subsection{Causal influence measures controllability}
\label{scontrol}

Causal influence is intimately related to control. Suppose an
experimenter wishes to understand interactions between components of a
complex system. For the causal DAG in Figure~\ref{DAGsobv}(d),
she is able to observe nodes $Y$ and $Z$, and manipulate node $X$. To
what extent can she control node $Y$? The notion of control has been
formalized information-theoretically in \cite{touchette04}:

%
\begin{Def}[(Perfect control)]
Node $Y$ is \emph{perfectly controllable} by node $X$ at $Z=z$ if,
given $z$,
\begin{longlist}[(ii)]
\item[(i)] states of $Y$ are a deterministic function of states of
$X$; and
\item[(ii)] manipulating $X$ gives rise to all states of $Y$.
\end{longlist}
\end{Def}

Perfect control can be elegantly characterized:

%
\begin{Thm}[(Information-theoretic characterization of perfect
controllability)]\label{tcontrol}
A node $Y$ with inputs $X$ and $Z$ is perfectly controllable by $X$
alone for $Z=z$ iff there exists a Markov transition matrix $R(x|z)$
such that
%
%
\renewcommand{\theequation}{C\arabic{equation}}
\setcounter{equation}{0}
\begin{eqnarray}
H\bigl(Y|z,{do} (x)\bigr):=\sum_x R(x|z) H\bigl(Y|z,{do} (x)\bigr)&=& 0\quad
\mbox{ and } 
\label{eC1}
\\
\sum_{x\in X} P\bigl(y|z, {do} (x)\bigr) R(x|z)&\neq&0\qquad
\mbox{for all }y. 
\label{eC2}
\end{eqnarray}
Here, $H(Y|z, {do} (x))$ denotes the conditional Shannon entropy of
$Y$, given that $Z=z$ has been observed and $X$ has been set to $x$.
\end{Thm}

\begin{pf}
The theorem restates the criteria in the definition. For a proof, see~\cite{touchette04}.
\end{pf}

It is instructive to compare Theorem~\ref{tcontrol} to our measure of
causal influence. The theorem highlights two fundamental properties of
perfect control. First, \eqref{eC1}, perfect control requires there is
no variation in $Y$'s behavior---aside from that due to the
manipulation via $X$---given that $z$ is observed. Second, \eqref
{eC2}, perfect control requires that all potential outputs of $Y$ can
be induced by manipulating node $X$. This suggests a measure of the
\emph{degree} of control\vadjust{\goodbreak} should reflect (i) the variability in $Y$'s
behavior that cannot be eliminated by imposing $X$ values and (ii) the
size of the repertoire of behaviors that can be induced on the target
by manipulating a source.

For the DAG under consideration, Theorem~\ref{thmdeco} states that
\[
\ci_{X\rightarrow Y}(P) = I(X;Y |Z)=H(Y|Z) - H(Y|X,Z).
\]
The first term, $H(Y|Z)$, quantifies size of the repertoire of outputs
of $Y$ averaged over manipulations of $X$. It corresponds to
requirement \eqref{eC2} in the characterization of perfect control:
that $P(y|z)>0$ for all $z$. Specifically, the causal influence,
interpreted as a measure of the degree of controllability, increases
with the size of the (weighted) repertoire of outputs that can be
induced by manipulations.

The second term, $H(Y|X,Z)$ [which coincides with $H(Y|Z, {do}(X))$
here], quantifies the variability in $Y$'s behavior that cannot be
eliminated by controlling $X$. It corresponds to requirement \eqref
{eC1} in the characterization of perfect control: that remaining
variability should be zero. Causal influence increases as the
variability $H(Y|Z, {do}(X))=\sum_z P(z) H(Y|z, {do}(X))$ tends toward
zero provided that the first term remains constant.
\end{appendix}

\section*{Acknowledgement}
We are grateful to G{\'a}bor Lugosi for a helpful hint for the proof of
Lemma~1 in Supplement S.1 and to Philipp Geiger for several corrections.

\begin{supplement}[id=suppA]
\stitle{Supplement to ``Quantifying causal influences''\\}
\slink[doi,text={10.1214/13-\break AOS1145SUPP}]{10.1214/13-AOS1145SUPP} 
\sdatatype{.pdf}
\sfilename{aos1145\_supp.pdf}
\sdescription{Three supplementary sections:\break (1)~Generating an i.i.d.
copy via random permutations; (2) Another option to define causal
strength; and (3) The problem of defining total influence.}
\end{supplement}

%

%


\printaddresses


\begin{thebibliography}{24}

\bibitem{Amari}
%
\begin{bbook}[auto:STB|2013/09/19|12:14:10]
\bauthor{\bsnm{Amari},~\bfnm{S.}\binits{S.}} \AND
\bauthor{\bsnm{Nagaoka},~\bfnm{H.}\binits{H.}}
(\byear{1993}).
\btitle{Methods of Information Geometry}.
\bpublisher{Oxford Univ. Press}, \blocation{New York}.
\bptok{imsref}%
\end{bbook}
%
\endbibitem

\bibitem{Avin}
%
\begin{bmisc}[auto:STB|2013/09/19|12:14:10]
\bauthor{\bsnm{Avin},~\bfnm{C.}\binits{C.}},
\bauthor{\bsnm{Shpitser},~\bfnm{I.}\binits{I.}} \AND
\bauthor{\bsnm{Pearl},~\bfnm{J.}\binits{J.}}
(\byear{2005}).
\bhowpublished{Identifiability of path-specific effects.
In \textit{Proceedings of the International Joint Conference in Artificial
Intelligence, Edinburgh, Scotland} 357--363. Professional Book Center, Denver}.
\bptok{imsref}%
\end{bmisc}
%
\endbibitem

\bibitem{AyKrakauer}
%
\begin{barticle}[auto:STB|2013/09/19|12:14:10]
\bauthor{\bsnm{Ay},~\bfnm{N.}\binits{N.}} \AND
\bauthor{\bsnm{Krakauer},~\bfnm{D.}\binits{D.}}
(\byear{2007}).
\btitle{Geometric robustness and biological networks}.
\bjournal{Theory in Biosciences}
\bvolume{125}
\bpages{93--121}.
\bptok{imsref}%
\end{barticle}
%
\endbibitem

\bibitem{AyInfoFlow}
%
\begin{barticle}[mr]
\bauthor{\bsnm{Ay},~\bfnm{Nihat}\binits{N.}} \AND
\bauthor{\bsnm{Polani},~\bfnm{Daniel}\binits{D.}}
(\byear{2008}).
\btitle{Information flows in causal networks}.
\bjournal{Adv. Complex Syst.}
\bvolume{11}
\bpages{17--41}.
\bid{doi={10.1142/S0219525908001465}, issn={0219-5259}, mr={2400125}}
\bptok{imsref}%
\end{barticle}
%
\endbibitem

\bibitem{cover}
%
\begin{bbook}[mr]
\bauthor{\bsnm{Cover},~\bfnm{Thomas~M.}\binits{T.~M.}} \AND
\bauthor{\bsnm{Thomas},~\bfnm{Joy~A.}\binits{J.~A.}}
(\byear{1991}).
\btitle{Elements of Information Theory}.
\bpublisher{Wiley}, \blocation{New York}.
\bid{doi={10.1002/0471200611}, mr={1122806}}
\bptok{imsref}%
\end{bbook}
%
\endbibitem

\bibitem{Granger1969}
%
\begin{barticle}[auto:STB|2013/09/19|12:14:10]
\bauthor{\bsnm{Granger},~\bfnm{C.~W.~J.}\binits{C.~W.~J.}}
(\byear{1969}).
\btitle{Investigating causal relations by econometric models and cross-spectral
methods}.
\bjournal{Econometrica}
\bvolume{37}
\bpages{424--38}.
\bptok{imsref}%
\end{barticle}
%
\endbibitem

\bibitem{Holland}
%
\begin{bincollection}[auto:STB|2013/09/19|12:14:10]
\bauthor{\bsnm{Holland},~\bfnm{P.~W.}\binits{P.~W.}}
(\byear{1988}).
\btitle{Causal inference, path analysis, and recursive structural equations
models}.
In \bbooktitle{Sociological Methodology}
(\beditor{\bfnm{C.}\binits{C.}~\bsnm{Clogg}}, ed.)
\bpages{449--484}.
\bpublisher{American Sociological Association}, \blocation
{Washington, DC}.
\bptok{imsref}%
\end{bincollection}
%
\endbibitem


\bibitem{Hoyer}
%
\begin{bincollection}[auto:STB|2013/09/19|12:14:10]
\bauthor{\bsnm{Hoyer},~\bfnm{P.}\binits{P.}},
\bauthor{\bsnm{Janzing},~\bfnm{D.}\binits{D.}},
\bauthor{\bsnm{Mooij},~\bfnm{J.}\binits{J.}},
\bauthor{\bsnm{Peters},~\bfnm{J.}\binits{J.}} \AND
\bauthor{\bsnm{Sch{\"o}lkopf},~\bfnm{B.}\binits{B.}}
(\byear{2009}).
\btitle{Nonlinear causal discovery with additive noise models}.
In \bbooktitle{Advances in Neural Information Processing Systems 21: 22nd Annual Conference on Neural Information Processing Systems 2008}
(\beditor{\bfnm{D.}\binits{D.}~\bsnm{Koller}},
\beditor{\bfnm{D.}\binits{D.}~\bsnm{Schuurmans}},
\beditor{\bfnm{Y.}\binits{Y.}~\bsnm{Bengio}} \AND
\beditor{\bfnm{L.}\binits{L.}~\bsnm{Bottou}}, eds.)
\bpages{689--696}.
\bpublisher{Curran Associates}, \blocation{Red Hook, NY}.
\bptok{imsref}%
\end{bincollection}
%
\endbibitem

\bibitem{suppcausalstrength}
%
\begin{bmisc}[auto:STB|2013/09/19|12:14:10]
\bauthor{\bsnm{Janzing},~\bfnm{D.}\binits{D.}},
\bauthor{\bsnm{Balduzzi},~\bfnm{D.}\binits{D.}},
\bauthor{\bsnm{Grosse-Wentrup},~\bfnm{M.}\binits{M.}} \AND
\bauthor{\bsnm{Sch{\"o}lkopf},~\bfnm{B.}\binits{B.}}
(\byear{2013}).
\bhowpublished{Supplement to ``Quantifying causal influences.''
DOI:\href{http://dx.doi.org/10.1214/13-AOS1145SUPP}{10.1214/}
\href{http://dx.doi.org/10.1214/13-AOS1145SUPP}{13-AOS1145SUPP}.}
\bptok{imsref}%
\end{bmisc}
%
\endbibitem

\bibitem{Lauritzen1996}
%
\begin{bbook}[mr]
\bauthor{\bsnm{Lauritzen},~\bfnm{Steffen~L.}\binits{S.~L.}}
(\byear{1996}).
\btitle{Graphical Models}.
\bseries{Oxford Statistical Science Series}
\bvolume{17}.
\bpublisher{Oxford Univ. Press}, \blocation{New York}.
\bid{mr={1419991}}
\bptok{imsref}%
\end{bbook}
%
\endbibitem

\bibitem{Lewontin}
%
\begin{barticle}[auto:STB|2013/09/19|12:14:10]
\bauthor{\bsnm{Lewontin},~\bfnm{R.~C.}\binits{R.~C.}}
(\byear{1974}).
\btitle{Annotation: The analysis of variance and the analysis of causes}.
\bjournal{American Journal Human Genetics}
\bvolume{26}
\bpages{400--411}.
\bptok{imsref}%
\end{barticle}
%
\endbibitem

\bibitem{Massey}
%
\begin{bmisc}[auto:STB|2013/09/19|12:14:10]
\bauthor{\bsnm{Massey},~\bfnm{J.}\binits{J.}}
(\byear{1990}).
\bhowpublished{Causality, feedback and directed information.
In \textit{Proc. 1990 Intl. Symp. on Info. Th. and Its Applications}.
Waikiki, Hawaii}.
\bptok{imsref}%
\end{bmisc}
%
\endbibitem

\bibitem{Northcott}
%
\begin{barticle}[auto:STB|2013/09/19|12:14:10]
\bauthor{\bsnm{Northcott},~\bfnm{R.}\binits{R.}}
(\byear{2008}).
\btitle{Can ANOVA measure causal strength?}
\bjournal{The Quaterly Review of Biology}
\bvolume{83}
\bpages{47--55}.
\bptok{imsref}%
\end{barticle}
%
\endbibitem

\bibitem{Pearl00}
%
\begin{bbook}[mr]
\bauthor{\bsnm{Pearl},~\bfnm{Judea}\binits{J.}}
(\byear{2000}).
\btitle{Causality: Models, Reasoning, and Inference}.
\bpublisher{Cambridge Univ. Press}, \blocation{Cambridge}.
\bid{mr={1744773}}
\bptok{imsref}%
\end{bbook}
%
\endbibitem

\bibitem{Pearlindirect}
%
\begin{bincollection}[auto:STB|2013/09/19|12:14:10]
\bauthor{\bsnm{Pearl},~\bfnm{J.}\binits{J.}}
(\byear{2001}).
\btitle{Direct and indirect effects}.
In \bbooktitle{Proceedings of the 17th Conference on Uncertainty in
Artificial Intelligence (UAI2001)}
\bpages{411--420}.
\bpublisher{Morgan Kaufmann}, \blocation{San Francisco, CA}.
\bptok{imsref}%
\end{bincollection}
%
\endbibitem


\bibitem{Perez-Cruz}
%
\begin{bmisc}[auto:STB|2013/09/19|12:14:10]
\bauthor{\bsnm{P{\'e}rez-Cruz},~\bfnm{F.}\binits{F.}}
(\byear{2009}).
\bhowpublished{Estimation of information theoretic measures for
continuous random
variables.
In \textit{Advances in Neural Information Processing Systems 21: 22nd Annual Conference on Neural
Information Processing Systems 2008} (D. Koller, D. Schuurmans, Y.~Bengio and
L. Bottou, eds.) 1257--1264. Curran Associates, Red Hook, NY.}
\bptok{imsref}%
\end{bmisc}
%
\endbibitem

\bibitem{Jonastpami}
%
\begin{barticle}[auto:STB|2013/09/19|12:14:10]
\bauthor{\bsnm{Peters},~\bfnm{J.}\binits{J.}},
\bauthor{\bsnm{Janzing},~\bfnm{D.}\binits{D.}} \AND
\bauthor{\bsnm{Sch{\"o}lkopf},~\bfnm{B.}\binits{B.}}
(\byear{2011}).
\btitle{Causal inference on discrete data using additive noise models}.
\bjournal{IEEE Transac. Patt. Analysis and Machine Int.}
\bvolume{33}
\bpages{2436--2450}.
\bptok{imsref}%
\end{barticle}
%
\endbibitem

\bibitem{UAIidentifiability}
%
\begin{bmisc}[auto:STB|2013/09/19|12:14:10]
\bauthor{\bsnm{Peters},~\bfnm{J.}\binits{J.}},
\bauthor{\bsnm{Mooij},~\bfnm{J.}\binits{J.}},
\bauthor{\bsnm{Janzing},~\bfnm{D.}\binits{D.}} \AND
\bauthor{\bsnm{Sch{\"o}lkopf},~\bfnm{B.}\binits{B.}}
(\byear{2001}).
\bhowpublished{Identifiability of causal graphs using functional
models. In
\emph{Proceedings of the 27th Conference on Uncertainty in Artificial
Intelligence (UAI 2011)} 589--598.
AUAI Press, Corvallis, OR.
Available at \url{http://uai.sis.pitt.edu/papers/11/p589-peters.pdf}}.
\bptok{imsref}%
\end{bmisc}
%
\endbibitem

\bibitem{Robins}
%
\begin{barticle}[auto:STB|2013/09/19|12:14:10]
\bauthor{\bsnm{Robins},~\bfnm{J.~M.}\binits{J.~M.}} \AND
\bauthor{\bsnm{Greenland},~\bfnm{S.}\binits{S.}}
(\byear{1992}).
\btitle{Identifiability and exchangeability for direct and indirect effects}.
\bjournal{Epidemiology}
\bvolume{3}
\bpages{143--155}.
\bptok{imsref}%
\end{barticle}
%
\endbibitem

\bibitem{Schreiber}
%
\begin{barticle}[auto:STB|2013/09/19|12:14:10]
\bauthor{\bsnm{Schreiber},~\bfnm{T.}\binits{T.}}
(\byear{2000}).
\btitle{Measuring information transfer}.
\bjournal{Phys. Rev. Lett.}
\bvolume{85}
\bpages{461--464}.
\bptok{imsref}%
\end{barticle}
%
\endbibitem

\bibitem{Spirtes}
%
\begin{bbook}[mr]
\bauthor{\bsnm{Spirtes},~\bfnm{Peter}\binits{P.}},
\bauthor{\bsnm{Glymour},~\bfnm{Clark}\binits{C.}} \AND
\bauthor{\bsnm{Scheines},~\bfnm{Richard}\binits{R.}}
(\byear{1993}).
\btitle{Causation, Prediction, and Search}.
\bseries{Lecture Notes in Statistics}
\bvolume{81}.
\bpublisher{Springer}, \blocation{New York}.
\bid{doi={10.1007/978-1-4612-2748-9}, mr={1227558}}
\bptok{imsref}%
\end{bbook}
%
\endbibitem

\bibitem{touchette04}
%
\begin{barticle}[mr]
\bauthor{\bsnm{Touchette},~\bfnm{Hugo}\binits{H.}} \AND
\bauthor{\bsnm{Lloyd},~\bfnm{Seth}\binits{S.}}
(\byear{2004}).
\btitle{Information-theoretic approach to the study of control systems}.
\bjournal{Phys. A}
\bvolume{331}
\bpages{140--172}.
\bid{doi={10.1016/j.physa.2003.09.007}, issn={0378-4371}, mr={2046360}}
\bptok{imsref}%
\end{barticle}
%
\endbibitem

\bibitem{ZhangUAI}
%
\begin{bmisc}[auto:STB|2013/09/19|12:14:10]
\bauthor{\bsnm{Zhang},~\bfnm{K.}\binits{K.}} \AND
\bauthor{\bsnm{Hyv{\"a}rinen},~\bfnm{A.}\binits{A.}}
(\byear{2009}).
\bhowpublished{On the identifiability of the post-nonlinear causal model.
In \textit{Proceedings of the 25th Conference on Uncertainty in Artificial
Intelligence, Montreal, Canada} 647--655. AUAI Press, Arlington, VA}.
\bptok{imsref}%
\end{bmisc}
%
\endbibitem

\end{thebibliography}
\end{document}